# Probability Bracket Notation: Probability Space, Conditional Expectation and Introductory Martingales

Dr. Xing M. Wang

## Table of Contents

## Abstract

In this paper, we continue to explore the consistence and usability of *Probability Bracket Notation (PBN)* proposed in our previous articles. After a brief review of *PBN* with dimensional analysis, we investigate probability spaces in terms of *PBN* by introducing probability spaces associated with random variables (*R.V*) or associated with stochastic processes (*S.P*). Next, we express several important properties of conditional expectation (*CE*) and some their proofs in *PBN*. Then, we introduce martingales based on sequence of *R.V* or based on filtration in *PBN*. In the process, we see *PBN* can be used to investigate some probability problems, which otherwise might need explicit usage of Measure theory. Whenever applicable, we use dimensional analysis to validate our formulas and use graphs for visualization of concepts in *PBN*. We hope this study shows that *PBN*, stimulated by and adapted from Dirac notation in Quantum Mechanics (QM), may have the potential to be a useful tool in probability modeling, at least for those who are already familiar with Dirac notation in QM.

## 1. Introduction

In our previous papers [1], we proposed *Probability Bracket Notation (PBN),* a new set of symbols for probability modeling, inspired by *Dirac notation*, which has been widely used in quantum mechanics (QM, [2]). We used *PBN* to formulas of basic probability theories, Markov chains and stochastic processes, represented the time evolution differential equation of time-continuous Markov chains in Schrodinger picture. Then, in Ref [3], we discussed the relationship between a sample space and its induced Hilbert space, investigated Diffusion maps based on Markov chains, and some phase spaces of many-body systems in Thermophysics.

In the above two articles, we never mentioned the concepts related to *probability space*, such as sigma-fields and filtration (see [4] or [5]). Now we want to investigate some advanced topics in probability theories and *S.P.* In this section, we give a brief review of *PBN* with the help of dimensional analysis, which shows that introducing Dirac delta function is the correct choice for continuous *R.V*, and it enable us to investigate a set of probability problems without explicitly using Measure theory (see [5], [7] or [9]).

### 1.1. A Brief Review of *PBN* and Dimensional Analysis

Before we start to discuss probability space, we first give a brief review of some basic *PBN* symbols, introduced in Ref. [1] for time-independent sample space. We assume *R.V.* *X* is an observable, the sample space Ω is the set of all possible outcomes of *X*.

1. ***Conditional Probability*** of event *A* under evidence *B* as a ***probability bracket*** from a ***P-bra*** and a ***P-ket*** :

$$P-bra: P(A| \qquad P-ket: |B) \qquad P-bracket: P(A|B) \equiv (A|B) \qquad (1.1.1a)$$

$$P(A|B) = 0 \ \ if \ \ A \cap B = \varnothing; \quad P(A|B) = 1 \ \ if \ \ A \supset B \neq \varnothing \qquad (1.1.1b)$$

$$P(A|B) = P(A|I|B), \quad \text{where } I \text{ is an unit operator} \qquad (1.1.1c)$$

2. ***The base***:

$$Discrete\ R.V.: \quad X|x_i) = x_i|x_i); \quad P(x_i|x_j) = \delta_{i,j} \qquad (1.1.2a)$$

$$Continuous\ R.V. \quad X|x) = x|x); \quad P(x|x') = \delta(x - x') \qquad (1.1.2b)$$

3. ***The unit operator*** (***Completeness***) :

$$Discrete\ R.V. \qquad I = \sum\nolimits_{x_i \in \Omega} |x_i)P(x_i| \qquad (1.1.3a)$$

$$Continuous\ R.V. \quad I = \int_{x \in \Omega} dx\,|x)P(x| \qquad (1.1.3b)$$

4. ***Probability of subset H*** in sample space Ω. In discrete and continuous cases:

$$P(H) \equiv P(H \mid \Omega) = \sum\nolimits_{x_i \in \Omega} (H \mid x_i) P(x_i \mid \Omega)$$
$$= \sum\nolimits_{x_i \in H} P(x_i \mid \Omega) \equiv \sum\nolimits_{x_i \in H} m(x_i) \qquad (1.1.4a)$$

$$P(H) \equiv P(H \mid \Omega) = \int_{x \in \Omega} dx (H \mid x) P(x \mid \Omega)$$
$$= \int_{x \in H} dx P(x \mid \Omega) \equiv \int_{x \in H} dx f(x) \qquad (1.1.4b)$$

5. ***Expectation value*** of $R.V.\ X$ (in discrete and continuous cases):

$$E[X] \equiv \langle X \rangle \equiv P(\Omega \mid X \mid \Omega) \underset{Eq.(1.1.3)}{=} \sum\nolimits_{x_i \in \Omega} P(\Omega \mid X \mid x_i) P(x_i \mid \Omega)$$
$$\underset{Eq.(1.1.2)}{=} \sum\nolimits_{x_i \in \Omega} P(\Omega \mid x_i) x_i P(x_i \mid \Omega) \underset{Eq.(1.1.1)}{=} \sum\nolimits_{x \in \Omega} x_i\, P(x_i \mid \Omega) \equiv \sum\nolimits_{x_i \in \Omega} x_i\, m(x_i) \qquad (1.1.5a)$$

$$E[X] \equiv \langle X \rangle \equiv P(\Omega \mid X \mid \Omega) \underset{Eq.(1.1.3)}{=} \int_{x \in \Omega} dx\, P(\Omega \mid X \mid x) P(x \mid \Omega)$$
$$\underset{Eq.(1.1.2)}{=} \int_{x \in \Omega} dx\, P(\Omega \mid x) x (x \mid \Omega) \underset{Eq.(1.1.1)}{=} \int_{x \in \Omega} dx\, x\, P(x \mid \Omega) \equiv \int_{x \in \Omega} dx\, x\, f(x) \qquad (1.1.5b)$$

6. ***Conditional expectation*** ($C.E$) value of $R.V.$ X under evidence $H$:

$$E[X \mid H] \equiv P(\Omega \mid X \mid H) \qquad (1.1.6)$$

7. If $A$ and $B$ are ***mutually independent***, then:

$$P(A \mid B) = P(A \mid \Omega) \qquad (1.1.7a)$$

And, if $X$ is independent of $H$, the $C.E$ of $X$ given $H$ becomes:

$$E[X \mid H] \equiv P(\Omega \mid X \mid H) = P(\Omega \mid X \mid \Omega) \equiv \langle X \rangle \qquad (1.1.7b)$$

***Comment 1.1.1:*** The definition of conditional probability is given by:

$$(A \mid B) \equiv P(A \mid B) = \frac{P(A \cap B)}{P(B)} = \frac{|A \cap B|}{|B|}$$

It requires that $P(B) > 0$, while for *continuous case*, from Eq. (1.1.4b) we might have $P(B) = 0$. But, assuming $|A| \equiv \int_{x \in \Omega} dx > 0$ and $B = x \in A \subseteq \Omega$, we can keep using Eq. (1.1.4b) and still get:

$$P(A \mid x) = \frac{P(A \cap x \mid \Omega)}{P(x \mid \Omega)} = \frac{P(x \mid \Omega)}{P(x \mid \Omega)} = \frac{f(x)}{f(x)} = 1 \qquad (1.1.8)$$

Now, if $A$ is also a point in continuous Ω, then we have a Dirac delta function $\delta(x - x')$ (see [2], §3.2) as in Eq (1.1.2b). It is consistent with our unit operator:

$$P(x|\Omega) \underset{(1.1.3b)}{=} \int_{x'\in\Omega} dx' P(x|x') P(x'|\Omega) \underset{(1.1.2b)}{=} \int_{x'\in\Omega} dx' \delta(x-x') P(x'|\Omega) = P(x|\Omega)$$

As for *discrete cases*, unless mentioned otherwise, we will always assume that $0 < |B|/|\Omega| \le 1$ as long as $B \neq \varnothing$.

***Dimensional Analysis***: In Physics, *dimensional analysis* plays a very important rule in verifying equations. We denote the *physical dimension* (*P.D)* of a quantity *Q* by $[Q]_{P.D}$. For examples, if *Q* is a position in one-dimensional physical space, then $[Q]_{P.D} = L$; if *Q* represents time, then $[Q]_{P.D} = T$. Assuming $[x]_{P.D} = L$ in our continuous cases (in most of our examples, *x* represents position), then from Eq. (1.1.3b) and (1.1.5b) we have:

***Continuous cases***:

$$[(x|]_{P.D} = [|x)]_{P.D} = [dx]^{-1/2} = L^{-1/2}, \quad [f(x)]_{P.D} = [dx]_{P.D}{}^{-1} = L^{-1} \tag{1.1.9}$$

Then we can find the *P.D* for other quantities as follows:

$$|\Omega) = \int_x dx|x)P(x|\Omega) = \int_x dx|x)f(x) \Rightarrow [|\Omega)]_{P.D} = L^{-1/2} \tag{1.1.10a}$$

$$P(\Omega| = \int_x dx\, P(\Omega|x)(x| = \int_x dx\,(x| \Rightarrow [(\Omega|]_{P.D} = L^{1/2}, \ \therefore [(\Omega|\Omega)]_{P.D} = 1 \tag{1.1.10b}$$

$$\begin{aligned} &\text{If } P(\Omega|A)_{P.D} = 1, \text{ then } [|A)]_{P.D} = L^{-1/2}; \\ &\text{If } [P(A|\Omega)]_{P.D} = [P(A)]_{P.D} = 1, \text{ then } [P(A|]_{P.D} = L^{1/2} \end{aligned} \tag{1.1.10c}$$

Had we written Eq. (1.1.2b) as $(x|x') = \delta_{x,x'}$, we would have: $(x|x) = 1 \Rightarrow [|x)]_{P.D} = 1 = L^0$ in contradiction to the *P.D* from Eq. (1.1.3b). Using dimensional analysis, we can combine Eq. (1.1.1), (1.1.7a), (1.1.8) as:

$$P(A|B) = 0, \quad \textit{if } A \cap B = \varnothing \tag{1.1.11a}$$

$$1'.\ P(A|B) \equiv (A|B): \ P(A|B) = 1, \quad \textit{if } A \supseteq B \supset \varnothing \textit{ and } [P(A|\Omega)]_{P.D} = 1 \tag{1.1.11b}$$

$$P(A|B) = P(A|\Omega), \quad \textit{if } A \textit{ and } B \textit{ are independent} \tag{1.1.11c}$$

In the phase space of Thermophysics (see [3]), we sometimes use another continuous variable *P* (*momentum*), having $[p]_{P.D} = MLT^{-1}$. For its base, we have:

$$[P(p|p')]_{P.D} = [\delta(p-p')]_{P.D} = [p]_{P.D}{}^{-1} = M^{-1}L^{-1}T \tag{1.1.12}$$

As for ***discrete cases***, we can easily check that following quantities are all *dimensionless*:

$$[|x_i)]_{P.D} = [P(x_i|]_{P.D} = [|\Omega)]_{P.D} = [P(\Omega|]_{P.D} = [m(x_i)]_{P.D} = 1 = L^0. \tag{1.1.13}$$

***Comment 1.1.2:*** With the help of Dirac delta function, *Bayes Formula* now is also valid for continuous base bracket defined in Eq. (1.1.2b):

$$P(x|x') = \frac{P(x'|x)P(x|\Omega)}{P(x'|\Omega)} = \frac{\delta(x'-x)P(x|\Omega)}{P(x'|\Omega)} = \delta(x'-x) \tag{1.1.14a}$$

$$\therefore P(x|x') \equiv (x'|x) = \delta(x-x'), \text{ as expected} \tag{1.1.14b}$$

## 2. Probability Space and *PBN*

In this section, we give a brief discussion of probability space associated with random variables or stochastic processes, as a reintroduction to *PBN*.

### 2.1. Basic Concepts: Probability Space and *R.V.*

The following statements and definitions 2.1.1-2.1.6 are quoted from §1.1-1.2, Ref. [5].

***Probability space*** $(\Omega, \Im, P)$ : The *sample space* $\Omega$ is a set of all possible outcomes $\omega \in \Omega$ of some random experiment. *Probabilities* are assigned by $A \mapsto P(A)$ to A in a subset $\Im$ of all possible sets of outcomes. The *event space* $\Im$ represents both the amount of information available as a result of the experiment conducted and the collection of all events of possible interest to us. A pleasant mathematical framework results by imposing on $\Im$ the structural conditions of *σ-field*.

We use $2^\Omega$ to denote the set of all possible subsets of Ω. The event space is thus a subset $\Im$ of $2^\Omega$, consisting of all allowed events, that is, those events to which we shall assign probabilities. We next define the structural conditions imposed on $\Im$.

**Definition 2.1.1**: We say that $\Im$ is a σ-field (or a σ-algebra), if
(a) $\Omega \in \Im$.
(b) If $A \in \Im$ then $A^c \in \Im$ as well (where $A^c \equiv \neg A \equiv \Omega - A$).
(c) If $A_i \in \Im$ for $i = 1, 2, \ldots$ then also $\bigcup_i A_i \in \Im$.

**Remark:** Using DeMorgan's law we can easily check that if $A_i \in \Im$ for $i = 1, 2, \ldots$ and $\Im$ is a σ-field, then also $\bigcap_i A_i \in \Im$.

**Definition 2.1.2**: A pair $(\Omega, \Im)$ with $\Im$ a σ-field of subsets of Ω is called a *measurable space*. Given a measurable space, a *probability measure P* is a function $P : \Im \mapsto [0,1]$, having the following properties:

1. $0 \le P(A)$ for all $A \in \Im$. In *PBN*, this is expressed as $0 \le P(A|\Omega)$ for all $A \in \Im$.
2. $P(\Omega) = 1$. In *PBN*, this is expressed as $P(\Omega|\Omega) = 1$.

3. $P(A)=\sum_{n=1}^{\infty}P(A_n)$, whenever $A=\bigcup_{n=1}^{\infty}A_n$ is a countable union of disjoint sets $A_n\in\Im$ (that is, $A_m\cap A_n=\varnothing$, for all $m\neq n$). In *PBN*, this is expressed as: $P(A\,|\,\Omega)=\sum_{n=1}^{\infty}P(A\,|A_n)P(A_n\,|\,\Omega)\underset{Eq.(1.1)}{=}\sum_{n=1}^{\infty}P(A_n\,|\,\Omega)$

**Definition 2.1.3:** A *probability space* is a *triplet* $(\Omega,\Im,P)$, with *P* a probability measure on the measurable space $(\Omega,\Im)$.

**Definition 2.1.4a**: A *Random Variable* (*R.V.*) is a function $X:\Omega\mapsto\Re$ such that for $\forall\alpha\in\Re$ the set $\{\omega: X(\omega)\leq\alpha\}$ is in $\Im$ (such a function is also called a *measurable function*).

**Definition 2.1.5**: We say that *R.V. X* and *Y* defined on the same probability space $(\Omega,\Im,P)$ are almost surely the same if $P(\{\omega: X(\omega)\neq Y(\omega)\})=0$. This shall be denoted by $X\overset{a.s.}{=}Y$. More generally, the same notation applies to any property of a *R.V.* For example, $X(\omega)\geq 0$ *a.s.* means that $P(\{\omega: X(\omega)<0\})=0$.

***Comment:*** Based on the above definitions, we can see the relationship between *probability space* and notations in *PBN* as follows:

1. $\Omega\mapsto P(\Omega\,|$. This is because $P(\Omega|$ represents all possible outcomes:
$$P(\Omega\,|=\sum_i P(\Omega\,|x_i)P(x_i\,|=\sum_i P(x_i\,|\equiv\sum_{x\in\Omega}P(x\,| \qquad (2.1.1a)$$
$$P(\Omega\,|=\int_{x\in\Omega}dx\,P(\Omega\,|\,x)P(x\,|=\int_{x\in\Omega}dx\,P(x\,| \qquad (2.1.1b)$$
2. $P: A\in\Im\mapsto P(A\,|\,\Omega)$. This is because we have:
$$P(A\,|\,\Omega)=\sum_i P(A\,|x_i)P(x_i\,|\,\Omega)=\sum_{\omega_i\in A}m(\omega_i)\equiv\sum_{x\in A}m(x)=P(A) \qquad (2.1.2a)$$
$$P(A\,|\,\Omega)=\int_{x\in A}dx\,P(A\,|\,x)P(x\,|\,\Omega)\equiv\int_{x\in A}dx\,f(x)=P(A) \qquad (2.1.2b)$$
3. $P:\Omega\in\Im\mapsto|\,\Omega)=\sum_i|\,x_i)P(x_i\,|\,\Omega)$. (2.1.3a)
$$P:\Omega\in\Im\mapsto|\,\Omega)=\int_x dx\,|\,x)P(x\,|\,\Omega) \qquad (2.1.3b)$$
4. Ω is all possible outcomes of observable *X* ↔The observable *X* is a R.V.
5. Observables *X* and *Y* have identical base P-kets ↔ $X\overset{a.s.}{=}Y$.

Next, let us present the definitions related to discrete-time martingales given by Ref. [4] (§9.5) and [5] (§4.1). Assume $(\Omega,\Im,P)$ is a probability space, **A** is a σ-field consisting of sub-sets of Ω, $X=X(\omega)$ is a function defined on Ω.

**Definition 2.1.4b**: *X* is called *measurable* on (Ω, **A** ) if for any real *x* we have:

$\{\omega, X(\omega)\leq x\}\in$ A

**Definition 2.1.6a**: The σ-field *generated* by *R.V. X* is the minimum σ-field on which *X* is measurable. It is denoted by $\Im(X)$.

**Definition 2.1.6b** ([4], definition 9.5.1; [5], definition 1.2.8): Given a *R.V. X* we denote by σ (*X*) the smallest σ-field $\mathcal{G} \subseteq \Im$ such that *X*(*ω*) is measurable on (Ω, $\mathcal{G}$ ). We call the σ-field generated by *X* and denote it σ (*X*) or $\Im_X \equiv \Im(X)$. Similarly, given *R.V.* $X_1, X_2, \dots X_n$ on the same measurable space $(\Omega, \Im)$, denote by $\sigma(X_k, k \le n)$ the smallest σ-field $\Im$ such that $X_k(\omega), k = 1, \dots n$ are measurable on $(\Omega, \Im)$. That is, $\sigma(X_k, k \le n)$ is the smallest σ-field containing $\sigma(X_k)$ for $k = 1, \dots n$.

**Lemma 2.1.1** (See [4], Lemma 9.5.2): If *X* is a discrete *R.V.*, its possible values are $x_1, x_2, \dots$, denote

$$A_i = \{\omega, X(\omega) = x_i\}, \quad i = 1, 2, \dots \tag{2.1.4}$$

Then $\Im(X)$ consists of $\varnothing, \Omega$ and all possible unions from elements of $\{A_i\}$

**Lemma 2.1.2** (See [4], Lemma 9.5.3): If *X* and *Y* are a discrete *R.V.*, their possible values are $x_1, x_2, \dots$ and $y_1, y_2, \dots$ denote

$$A_{ij} = \{\omega, X(\omega) = x_i, Y(\omega) = y_j\}, \quad i, j = 1, 2, \dots \tag{2.1.5}$$

Then $\Im(X, Y)$ consists of $\varnothing, \Omega$ and all possible unions from elements of $\{A_{ij}\}$

***Comment***: In PBN, this can be understood as:

$$\Im(X) \mapsto |\Omega_x) \text{ (one } R.V.\text{)}; \quad \Im(X, Y) \mapsto |\Omega_{x,y}) \text{ (two } R.V.\text{)} \tag{2.1.6}$$

**Definition 2.1.7 *Independence*:** We say that two events $A, B \in \Im$ are *P*-mutually independent if $P(A \cap B) = P(A)P(B)$, or $P(A \cap B | \Omega) = P(A | \Omega)P(B | \Omega)$.

In *PBN*, if two events $A, B \in \Im$ are *P*-mutually independent, then it implies that the original sample space is a direct product of two sample spaces such that:

$$\Omega = \Omega_1 \otimes \Omega_2, \; A \in \Omega_1, \; B \in \Omega_2, \; |\Omega) = |\Omega_1)|\Omega_2) \tag{2.1.7}$$

$$P(A \cap B | \Omega) = P([A \otimes \Omega_2] \cap [\Omega_1 \otimes B] | \Omega) = P(A \otimes B | \Omega) = P(A | \Omega)P(B | \Omega) \tag{2.1.8}$$

**Definition 2.1.8 *Borel (measurable) function*:** A function $g : \Re \mapsto \Re$ is called Borel (measurable) function if g is a *R.V.* on $(\Re, \mathcal{B})$. One can show that every continuous function *g* is Borel measurable. Further, every piecewise constant function *g* is also Borel measurable (where g piecewise constant means that *g* has at most countably many jump points between which it is constant).

**Definition 2.1.9 *Indicator function***: For any $A \in \Im$, an *indicator function* of *A* is defined:

$$\text{If } [P(A|\Omega)]_{P.D} = 1: \quad \boldsymbol{1}_A(\omega) = \begin{cases} = 1, \text{ if } \omega \in \text{A} \\ = 0, \text{ if } \omega \notin \text{A} \end{cases}; \quad \text{otherwise}: \quad \boldsymbol{1}_x(x') = \delta(x - x') \qquad (2.1.9a)$$

Indicator function $\boldsymbol{1}_A(\omega)$ is a *R.V.*, since

$$\{\omega : \boldsymbol{1}_A(\omega) \le x\} = \begin{cases} \Omega, \text{ if } x \ge 1 \\ \overline{A}, \text{ if } 0 \le x < 1 \\ \varnothing, \text{ if } x < 0 \end{cases} \Rightarrow \quad \text{For } \forall x \in \Re, \ \{\omega : \boldsymbol{1}_A(\omega) \le x\} \in \Im \qquad (2.1.9b)$$

Using *PBN*, we can express the indicator functions of $A$ as a projection operator $\boldsymbol{1}_A$ :

$$\begin{aligned} &\textit{Discrete R.V.} \quad \boldsymbol{1}_A = \sum\nolimits_{x_i \in A} | x_i ) P(x_i | \\ &\textit{Continuous R.V.} \ \ \boldsymbol{1}_A = \int_{x \in A} dx \, | x ) P(x | \end{aligned} \qquad (2.1.10)$$

We can transform the probability of set $A$ to the expectation value of operator $\boldsymbol{1}_A$ by using the following equivalence (see Eq. (3.2.5) of [4]):

$$P(A|\Omega) = P(\Omega | \boldsymbol{1}_A | \Omega) \equiv \langle \boldsymbol{1}_A \rangle \qquad (2.1.11)$$

This can be proved easily by using *PBN*. For example, in continuous case, we have:

$$\langle \boldsymbol{1}_A \rangle = P(\Omega | \boldsymbol{1}_A | \Omega) \underset{(2.1.10)}{=} \int_{x \in A} dx \, P(\Omega | x) P(x | \Omega) \underset{(1.1.1b)}{=} \int_{x \in A} dx \, P(x | \Omega) = P(A | \Omega) \qquad (2.1.12)$$

Note, Eq. (2.2.12) is true even if $A$ has only one point $x'$: $\langle \boldsymbol{1}_{x'} \rangle = P(x' | \Omega) \equiv f(x')$ .

## 2.2. Probability Space Associated with Random Variables

From section §2.1, we see that in current probability theories, we have the following definition chains:

$$\Omega \Rightarrow (\Omega, \Im) \Rightarrow P : A \in \Im \mapsto P(A) \in \Re \Rightarrow (\Omega, \Im, P) \Rightarrow X \ is \ R.V. on \, (\Omega, \Im, P)$$

But in *PBN*, we actually have proposed a quite different definition chain:

$$Observable \, X \Rightarrow base : \{| x)\} \Rightarrow \Omega \Rightarrow (\Omega, \Im) \Rightarrow P : x \in \Omega \mapsto P(x | \Omega) \in \Re \Rightarrow (\Omega, \Im, P)$$

Because of the way we define the probability space, an observable is automatically a *R.V.* on the probability space, i.e., an observable is equivalent to a *R.V.*

**One Observable Case:** Let *X* be an observable (or a $R.V.$). A probability space $(\Omega, \Im, P)$ is *associated* with *X*, if we have following properties:

**1. Observable $X \rightarrow \Omega$:** Ω is the space of all possible outcomes of *X*. For each element *x* (an outcome, or a base event) in the sample space Ω, *X* has a definite value *x*. We use a set of *base p-kets* to represent this property:

$$X|x) \equiv X|X=x) = x|x), \quad P(x|X = P(x|x \tag{2.2.1a}$$
$$x \in \Omega \Leftrightarrow \quad X|x) = x|x) \text{ and } P(\Omega|x) = 1 \tag{2.2.1b}$$
$$x \in \Omega \text{ and } g \text{ is a Borel function} \Leftrightarrow \quad g(X)|x) = g(x)|x) \tag{2.2.1c}$$

**2. Orthonormality of base:** All these outcomes are mutually disjoint ($x \cap x' = \emptyset$, if $x \neq x'$), so the conditional probability of any two outcomes $x, x' \in \Omega$ satisfies the following orthonormal relations as in Eq. (1.1.3):

$$\begin{aligned} &P(x|x') = \delta(x-x'), \quad P(x|X|x') = x\delta(x,x') \quad \text{(continuas } R.V.) \\ &P(x|x') = \delta_{x,x'} \quad P(x|X|x') = x\delta_{x,x'} \qquad \text{(discrete } R.V.) \end{aligned} \tag{2.2.2}$$

**3. Probability Existence**: The probability P of the probability space is defined by the probability distribution (or density) for a given outcome $x \in \Omega$ exists and is given by:

$$\begin{aligned} &P(x|\Omega) \equiv f(x) \geq 0, \quad \int_{x\in\Omega} dx\, P(x|\Omega) = 1 \quad \text{(continuas } R.V.) \\ &P(x|\Omega) \equiv m(x) \geq 0, \quad \sum_{x\in\Omega} P(x|\Omega) = 1 \qquad \text{(discrete } R.V.) \end{aligned} \tag{2.2.3a}$$

The probability of any set *H* in the *σ-field* $\Im$ on Ω can be obtained from the probability distribution (or density) by:

$$\begin{aligned} &P(H|\Omega) \equiv P(H) = \int_{x\in H} dx\, P(x|\Omega) \quad \text{(continuas } R.V.) \\ &P(H|\Omega) = P(H) = \sum_{x\in H} P(x|\Omega) \qquad \text{(discrete } R.V.) \end{aligned} \tag{2.2.3b}$$

**4. Completeness of Base**: Eq. (2.2.3) implies the completeness of base p-kets:

$$\begin{aligned} &\int_{x\in\Omega} |x)\, dx\, P(x| = 1 \qquad \text{(continuous } R.V.) \\ &\sum_{x\in\Omega} |x)\, P(x| = 1 \qquad \text{(discrete } R.V.) \end{aligned} \tag{2.2.4}$$

**5. Expectation value**: The expectation value of a Borel function (Definition 2.1.8) of g(*x*) now can be obtained by:

$$E[g(X)] \equiv \langle g(X)\rangle = P(\Omega\,|\,g(X)\,|\,\Omega) \underset{(2.2.4)}{=} \int_{x\in\Omega} dx\, P(\Omega\,|\,g(X)\,|\,x)P(x\,|\,\Omega)$$
$$\underset{(2.2.1)}{=} \int_{x\in\Omega} dx\, g(x)P(x\,|\,\Omega) \qquad \text{(continuas } R.V.\text{)} \tag{2.2.5a}$$

$$E[g(X)] \equiv \langle g(X)\rangle = P(\Omega\,|\,g(X)\,|\,\Omega) \underset{(2.2.4)}{=} \sum_{x\in\Omega} P(\Omega\,|\,g(X)\,|\,x)P(x\,|\,\Omega)$$
$$\underset{(2.2.1)}{=} \sum_{x\in\Omega} g(x)P(x\,|\,\Omega) \qquad \text{(discrete } R.V.\text{)} \tag{2.2.5b}$$

In particularly, the ***characteristic function*** of a distribution function (density) can be represented as:

$$\varphi(k) = \langle \mathrm{e}^{ikX}\rangle = P(\Omega\,|\,\mathrm{e}^{ikX}\,|\,\Omega) = \sum\nolimits_{x\in\Omega} \mathrm{e}^{ikx}\, P(x\,|\,\Omega) \quad \text{(discrete } R.V.\text{)} \tag{2.2.5c}$$
$$\varphi(k) = \langle \mathrm{e}^{ikX}\rangle = P(\Omega\,|\,\mathrm{e}^{ikX}\,|\,\Omega) = \int_{x\in\Omega} dx\, \mathrm{e}^{ikx}\, P(x\,|\,\Omega) \quad \text{(continuous } R.V.\text{)} \tag{2.2.5d}$$

From the above discussion, we can conclude that, for a probability space $(\Omega, \Im, P)$ associated with random variable $X$, the sample space $\Omega$ is all outcomes of measuring $X$, and probability $P$ for $\Omega$ is given by the distribution or density $P(x\,|\,\Omega)$ for each element $x\in\Omega$, and any set $H\in\Im$ will have a probability given by $P(H\,|\,\Omega)$. Therefore, all the information about the probability space is contained in the system p-ket $|\,\Omega)$. We can have the following statements for continuous $R.V$:

$$\Omega \equiv \text{ Range of } X, \text{ or all posibble outcomes of } X, \tag{2.2.6a}$$
$$P_X : \Omega \mapsto |\,\Omega) \Rightarrow \text{ for } \forall x\in\Omega,\, P(x) = P(x\,|\,\Omega), \tag{2.2.6b}$$
$$\text{For } \forall H\in\Im,\, P(H) = P(H\,|\,\Omega) = \int_{x\in H} dx P(x\,|\,\Omega) \tag{2.2.6c}$$

This is similar to the system state-ket of Hilbert space in 1D QM problem for $x\in\Re$:

$$\hat{x}\,|\,x\rangle = x\,|\,x\rangle, \quad \langle x\,|\,x'\rangle = \delta(x-x'), \quad \int_{-\infty}^{\infty} dx\,|\,x\rangle\langle x\,| = 1, \tag{2.2.7a}$$
$$P : \Re \mapsto |\,\Psi\rangle \Rightarrow \text{ for } \forall x\in\Re, P(x) = |\langle x\,|\,\Psi\rangle|^2 = |\,\Psi(x)\,|^2, \tag{2.2.7b}$$
$$\text{for } \forall A = [a,b]\in\Re, P(A) = \int_a^b dx\,|\langle x\,|\,\Psi\rangle|^2 = \int_a^b dx\,|\,\Psi(x)\,|^2 \tag{2.2.7c}$$

The dimensional analogous can be derived from (2.2.6) and (2.2.7):

$$\text{Assuming } [x]_{P.D} = L : [|\,x\rangle]_{P.D} = [|\,x)]_{P.D} = L^{-1/2}, \quad [P(x\,|\,\Omega)]_{P.D} = [|\langle x\,|\,\Psi\rangle|^2]_{P.D} = L^{-1}$$

Based on this similarity, we have discussed induced Hilbert space and induced sample (or probability) space in Ref. [2].

**Two Continuous Observables Case:** If $X$ and $Y$ are two continuous observables and a probability space $(\Omega, \Im, P)$ is associated with them, then there exists a joint density: for any joint event $(x \cap y) \in \Im$, we have a probability density $f(x, y)$ in *probability measure P* such that:

$$\begin{aligned} &|x,y) \equiv |x \cap y) \equiv |X = x \cap Y = y), \quad X|x,y) = x|x,y), \quad Y|x,y) = y|x,y) \\ &P(\Omega|x,y) = 1, \quad |\Omega) \equiv |\Omega_{x,y}), \quad P(\Omega| = (\Omega_{x,y}|, \\ &P(x,y|x',y') = \delta(x-x')\delta(y-y'), \quad P(x,y|\Omega) \equiv f(x,y) \geq 0 \end{aligned} \quad (2.2.8a)$$

$$\int |x,y)\, dxdy P(x,y| = 1, \quad P(\Omega|\Omega) = \int dxdy\, P(x,y|\Omega) = \int dxdy\, f(x,y) = 1 \quad (2.2.8b)$$

$$P(x|\Omega) \equiv P(x,*|\Omega) = \int dy\, P(x,y|\Omega), \quad P(y|\Omega) \equiv P(y,*|\Omega) = \int dx\, P(x,y|\Omega) \quad (2.2.8c)$$

Based on ***Comment*** 1.1.1, we have the following physical dimensions for two continuous observables:

$$[P(x,y|]_{P.D} = [|x,y)]_{P.D} = [dx\, dy]^{-1/2} = L^{-1}, \quad [f(x,y)]_{P.D} = [dx\, dy]_{P.D}{}^{-1} = L^{-2}. \quad (2.2.8d)$$

Because the base is complete, if $H$ is any set in $\Im$, we have:

$$\begin{aligned} &XY|H) = \int XY|x,y)\, dxdy\, P(x,y|H) = \int xy|x,y) dxdy P(x,y|H) = YX|H) \\ &\therefore (XY - YX)|H) \equiv [X,Y]|H) = 0, \quad \textit{for any } H \in \Im \end{aligned} \quad (2.2.9a)$$

Here we have used the definition of *commutator* in QM (see [2], §5.2). As we proposed in Ref [2], this means that, the existence of joint density means that, in the induced Hilbert space, the corresponding operators must have a complete common eigenvectors, and therefore they must commute each other:

$$\begin{aligned} &\hat{X}|x,y\rangle = x|x,y\rangle, \ \hat{Y}|x,y\rangle = y|x,y\rangle \\ &[\hat{X} - \hat{Y}] \equiv \hat{X}\hat{Y} - \hat{Y}\hat{X} = 0 \end{aligned} \quad (2.2.9b)$$

Similarly to Eq. (2.2.6), we can write:

$$\Omega \equiv \Omega_{x,y} = \text{Range of } (X,Y), \text{ or all posibble outcomes of } (X,Y), \quad (2.2.10a)$$

$$\begin{aligned} &P_{X,Y} : \Omega \mapsto |\Omega) \Rightarrow \text{ for } \forall x,y \in \Omega, P(x,y) = P(x,y|\Omega) = f(x,y), \\ &\text{for } \forall H \in \Im, P(H) = P(H|\Omega) = \int_{x,y \in H} dx\, dy P(x,y|\Omega) \end{aligned} \quad (2.2.10b)$$

This is analogue to the system state-ket of Hilbert space in 2D QM problem with $\vec{r} = (x,y)$:

$$\hat{x}_i|\vec{r}\rangle = x_i|\vec{r}\rangle, \quad \langle\vec{r}|\vec{r}\,'\rangle = \delta(\vec{r} - \vec{r}\,'), \quad \int_{-\infty}^{\infty}\int_{-\infty}^{\infty} d^2\vec{r}\,|\vec{r}\rangle\langle\vec{r}| = 1 \quad (2.2.11a)$$

$$P:\Re^2 \mapsto |\Psi\rangle \Rightarrow \text{ for } \forall(x,y)\in\Re^2, P(x,y) = |\langle x,y|\Psi\rangle|^2 = |\Psi(x,y)|^2, \quad (2.2.11b)$$

$$\text{for } \forall A = \{x\in[a,b], y\in[c,d]\}\in\Re^2, P(A) = \int_a^b dx\int_c^d dy\,|\Psi(x,y)|^2 \quad (2.2.11b)$$

If $X$ and $Y$ are *two independent* observables on a probability space $(\Omega,\Im,P)$, then there exist two subsets $\Omega_x$ and $\Omega_y$ in $\Im$ such that $\Omega = \Omega_x \otimes \Omega_y$, $\Omega_x \cap \Omega_y = \varnothing$, or:

$$|\Omega) = |\Omega_{x,y}) = |\Omega_x)\,|\Omega_y), \quad P(\Omega| = P(\Omega_{x,y}| = P(\Omega_x|\,P(\Omega_y| \quad (2.2.12a)$$

For any joint event $(x,y) \equiv (X = x \wedge Y = y) \equiv (X = x \otimes Y = y) \equiv (x \otimes y) \in \Im$, we have:

$$|x,y) = |x)\,|y), \quad P(x,y| = P(x|\,P(y| \quad (2.2.12b)$$

$$P(x|\Omega_x) \equiv f_X(x), \quad P(y|\Omega_y) \equiv f_Y(x) \quad (2.2.12c)$$

$$P(\Omega|x,y) = 1, \quad P(x,y|\Omega) = P(x|\Omega_x)P(y|\Omega_y) = f_X(x)\,f_Y(y) \quad (2.2.12d)$$

Such interpretation is consistent with *dimensional analysis*:

$$[|x,y)]_{P.D} = [|x)]_{P.D}[|y)]_{P.D} = L^{-1}; \quad (2.2.13a)$$

$$[|\Omega_{x,y})]_{P.D} = [|\Omega_x)\,|\Omega_y)]_{P.D} = L; \quad [f(x,y)]_{P.D} = [f_X(x)\,f_Y(y)]_{P.D} = L^2 \quad (2.2.13b)$$

In Hilbert space, we have the similar representations. For example, a system of two non-interacting spin-zero particles has the following state ket:

$$\langle\vec r_1,\vec r_2|\Psi\rangle = \langle\vec r_1|\Psi_1\rangle\langle\vec r_2|\Psi_2\rangle = \Psi_1(\vec r_1)\Psi_2(\vec r_2) \quad (2.2.14)$$

In general, if we have $n$ observables, $\{X_i, i = 1,2,...n\}$, then the associated probability space has following properties:

$$|\Omega) \equiv |\Omega_{x_1,...x_n}), P(\Omega| = P(\Omega_{x_1,...x_n}|, P(x_1,...,x_n|X_j|x_1',...,x_n') = x_j\prod_{i=1}^n \delta(x_i - x_i') \quad (2.2.15a)$$

$$P(x_1|\Omega) = \int dx_2\, P(x_1,x_2|\Omega) = \cdots = \int dx_2 dx_3\ldots dx_n P(x_1,x_2,...,x_n|\Omega) \quad (2.2.15b)$$

If the $n$ observables are mutually independent, then we have:

$$|\Omega) \equiv \prod_{i=1}^n |\Omega_{x_i}), \quad P(x_1,...,x_n|\Omega) = \prod_{i=1}^n P(x_i|\Omega_{x_i}) \quad (2.2.16)$$

## 2.3. Basic Concepts: Probability Space and Stochastic Process

Stochastic processes (*S.P.*) can be defined via *R.V.*-s ([4-5]) or via *filtered probability space* ([5-6]). We start with the following definition ([4], §5.1.1 and [5], §3.1) of *S.P.* based on *R.V.*-s.

**Definition 2.3.1**: Given $(\Omega, \Im, P)$, a *stochastic process* (*S.P.*) $\{X_t\}$ is a collection $\{X_t : t \in \mathrm{I}\}$ of *R.V.*-s where the index t belongs to the index set $\mathcal{I}$. Typically, $\mathcal{I}$ is an interval in $\Re$ (in which case we say that $\{X_t\}$ is a continuous time stochastic process), or a subset of $\{1, 2, \ldots n, \ldots\}$ (in which case we say that $\{X_t\}$ is a discrete time stochastic process. For a fixed $t \in \mathrm{I}$, we have an ordinary *R.V.*, while for a fixed ω, we call $t \mapsto X_t(\omega) \equiv X(\omega, t)$ the sample function (or sample path) of the *S.P.*

Recall our notation $\sigma(X_t)$ for the σ-field generated by $X_t$. The discrete time stochastic processes are merely countable collections of *R.V.*-s $X_1, X_2, X_3, \ldots$ defined on the same probability space. All relevant information about such a process during a finite time interval $\{1, 2, \ldots, n\}$ is conveyed by the σ-field $\sigma(X_1, X_2, \ldots, X_n)$, namely, the σ-field generated by the "rectangle" sets $\bigcap_{i=1}^{n} \{\omega : X_i(\omega) \le \alpha\}$ for $\alpha \in \Re$ (compare with *Definition 2.1.6b* of $\sigma(X)$). To deal with the full infinite time horizon we just take the σ-field $\sigma(X_1, X_2, \ldots)$ generated by the union of these sets over $n = 1, 2, \ldots$ It is not hard to verify that in this setting the σ-field $\sigma(X_1, X_2, \ldots)$ coincides with the smallest σ-field containing $\sigma(X_t)$) for all $t \in \mathrm{I}$, which we denote hereafter by $\Im_X \equiv \Im(X)$.

*S.P.* can also be defined through *filtered probability space*. A *filtration* represents any procedure of collecting more and more information as time goes on.

**Definition 2.3.1**: A ***filtration*** is a non-decreasing family of sub-σ-fields $\{\Im_n\}$ of our measurable space $(\Omega, \Im)$. That is, $\Im_0 \subseteq \Im_1 \subseteq \Im_2 \subseteq \cdots \Im_n \subseteq \cdots \subseteq \Im$ and $\Im_n$ is a σ-field for each *n*.

Given a filtration, we are interested in *S.P.* such that for each *n* the information gathered by that time suffices for evaluating the value of the *n*-th element of the process.

**Definition 2.3.2**: A *S.P.* $\{X_n, n = 0, 1, \ldots\}$ is ***adapted*** to a filtration $\{\Im_n\}$ if $\omega \mapsto X_n(\omega)$ is a *R.V.* on $(\Omega, \Im)$ for each *n*, that is, if $\sigma(X_n) \subseteq \Im_n$ for each *n*.

At this point you should convince yourself that $\{X_n\}$ is adapted to the filtration $\{\Im_n\}$ if and only if $\sigma(X_0, X_1, \ldots, X_n) \subseteq \Im_n$ for all n. That is,

**Definition 2.3.3**: The filtration $\{\mathcal{G}_n\}$ with $\mathcal{G}_n = \sigma(X_0, X_1, \ldots, X_n)$ is the minimal filtration with respect to which $\{X_n\}$ is adapted. We therefore call it the *canonical filtration* for the *S.P.* $\{X_n\}$.

Whenever clear from the context what it means, we shall use the notation $X_n$ both for the whole *S.P.* $\{X_n\}$ and for the *n*-th *R.V.* of this process, and likewise we may sometimes use $\Im_n$ to denote the whole filtration $\{\Im_n\}$.

We can interpolate $\Omega_n$ in *PBN* as the snapshot of the sample space *at* time *n*, and $\Im_n$ as the set of all snapshots taken *up to* time *n*. A measurement of *X* at time *n* is taken a snapshot of sample space at time *t* as shown in *Fig. 2.3.1.* And we can represent discrete filtration in *PBN* as:

$$\Im(X_n) \mapsto |\,\Omega_n), \quad \Im_n \mapsto |\,\Omega_0,...\Omega_n), \quad \Im \mapsto |\,\Omega) = |\,\Omega_0,...,\Omega_n,...) \tag{2.3.1}$$

Next, let us introduce the definition for continuous-time *S.P.*, based on filtration ([6], §1.1).

To define a stochastic process we need to add a time index *t*. This index can be discrete or continuous, finite or infinite. In what follows we will be dealing (almost) exclusively with continuous-time processes that start at date 0 and have an infinite horizon, so we will focus here on the case $t \in [0,\infty)$. We keep $\Omega, \Im$ and *P* as before, but in addition we add an increasing family of σ-algebras $\mathsf{F} = \{\Im_t, t \geq 0\}$ contained in $\Im$. That is,

$$\Im_s \subseteq \Im_t, \text{ all } s \leq t, \quad \text{and} \quad \Im_t \subseteq \Im, \text{ all } t \geq 0 \tag{2.3.2}$$

where $\Im = \Im_\infty$ is the smallest σ-algebra containing all the $\Im_t$'s. The family is called a *filtration*, and $(\Omega, \mathsf{F}, P)$ is called a *filtered probability space*. The interpretation is that $\Im_t$ is the set of events known at time t.

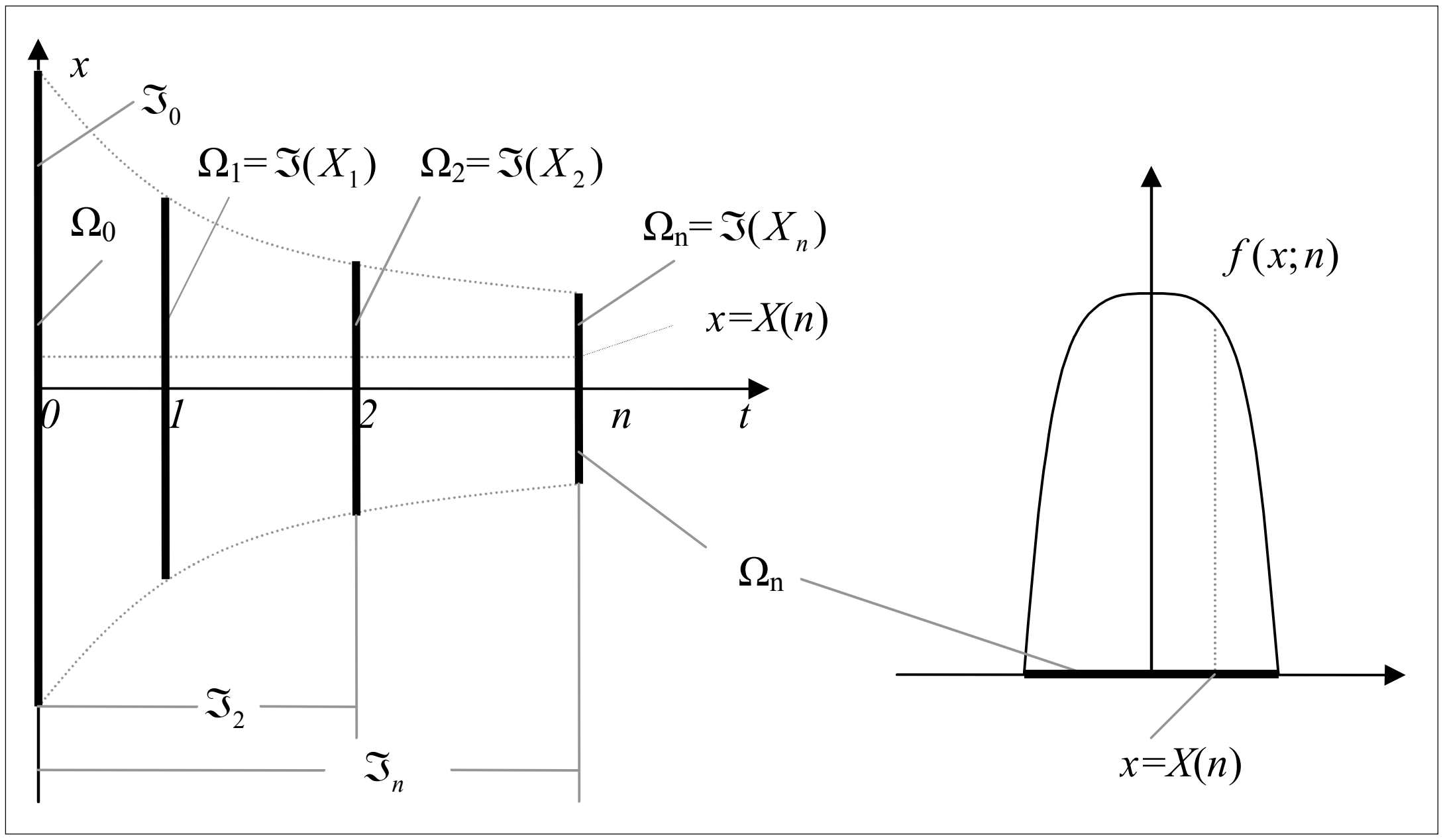

Fig. 2.3.1: Discrete $\Omega_n, \Im_n$ and $(X_n = x\,|\,\Omega) = (x\,|\,\Omega_n) = f(x;n)$

A *stochastic process* is a function on a filtered probability space with certain measurability properties. Specifically, let $(\Omega, \mathsf{F}, P)$ be a filtered probability space with time index $t \in [0,\infty) = \Re_+$, and let $\mathcal{B}_+$ denote the Borel subsets of $\Re_+$. A continuous-time stochastic process is a mapping $X:[0,\infty)\times\Omega \mapsto \Re$ that is measurable with respect to $\mathcal{B}_+ \times \Im$. That is, *X* is jointly measurable in (*t*, ω). Given a probability space $(\Omega, \Im, P)$ and a filtration $\mathsf{F} = \{\Im_t, t \ge 0\}$, we say that the stochastic process $X:[0,\infty)\times\Omega \mapsto \Re$ is *adapted* to $\mathfrak{F}$ if $X(\omega,t)$ is $\Im_t$-measurable, all t.

For each fixed $t \in [0,\infty) = \Re_+$, the mapping $X(\omega,.):[0,\infty) \mapsto \Re$ is an ordinary random variable on the probability space $(\Omega, \Im_t, P_t)$, where $P_t$ is the restriction of *P* to $\Im_t$. That is, $X(\omega,.)$ is an $\Im_t$-measurable function of ω. For each fixed $\omega \in \Omega$, the mapping $X:[0,\infty)\times\Omega \mapsto \Re$ is a Borel-measurable function of *t*. The mapping $X(.,t)$ is called a *realization* or *trajectory* or *sample path*.

From the above statements, we see that the standard definition chain is:

$$\Omega \Rightarrow (\Omega,\Im_t) \Rightarrow P_t : A \in \Im_t \mapsto P_t(A) \in \Re \Rightarrow (\Omega,\Im_t,P_t) \Rightarrow X_t \textit{ is R.V. on } (\Omega,\Im_t,P_t)$$
$$\Rightarrow (\Omega,\Im_t,P_t) \Rightarrow (\Omega,\mathsf{F},P) \Rightarrow X(\omega,t) \textit{ is adapted to } \mathsf{F} \qquad (2.3.3)$$

In *PBN*, we actually have a quite different definition chain:

$$\textit{Observable } X(t)\ (t \in \mathsf{I}\ ) \Rightarrow \textit{base}: \{|\, x;t)\} \Rightarrow \Omega_t \Rightarrow P_t, x \in \Omega_t \mapsto P(x \,|\, \Omega_t) \equiv f(x,t) \in \Re$$
$$\Rightarrow (\Omega_t,\Im_t,P_t) \Rightarrow P_t : A \subseteq \Omega_t \mapsto P(A \,|\, \Omega_t) \Rightarrow (\Omega,\mathsf{F},P) \Rightarrow X(\omega,t) \textit{ is adapted to } \mathsf{F} \qquad (2.3.4)$$

The notations for Ω, *X*(*t*) and $\Omega_t$ are as follows:

(a). For any $H \in \Im_t$, $(\Omega \,|\, H) = 1$; $P(\Omega \,|\, \Omega_t) = 1 \because \Omega_t \subset \Im_t$; $X(t)\,|\,\Omega) = X(t)\,|\,\Omega_t)$
(b). For discrete time: $|\,\Omega) = |\,\Omega_0,\Omega_1,\Omega_2,...) = |\,\Omega_{0,1,2,...})$, $P(x;n \,|\, \Omega) = P(x \,|\, \Omega_n) \equiv f(x;n)$
(c). For continuous time $t \in \mathsf{I} \subseteq [0,\infty)$: $|\,\Omega) = |\,\Omega_{\mathsf{I}}\ )$; $P(x;t \,|\, \Omega) = P(x \,|\, \Omega_t) \equiv f(x;t)$
(d). Ω is the collection of all possible outcomes of *X(t)* of all time $t \in \mathsf{I} \subseteq \Re_+$. When time changes, this collection do not change; only the probability distribution changes. This means that for $\forall t \in \mathsf{I} \subseteq \Re_+$ if $x \in \Omega_t$, then $x \in \Omega$.
(e). We define $\Omega_t$ as a subspace with a base consisting of all possible outcomes of *X(t)* at given time t, therefore $\Omega_t \subseteq \Im(X_t) \subset \Im$.

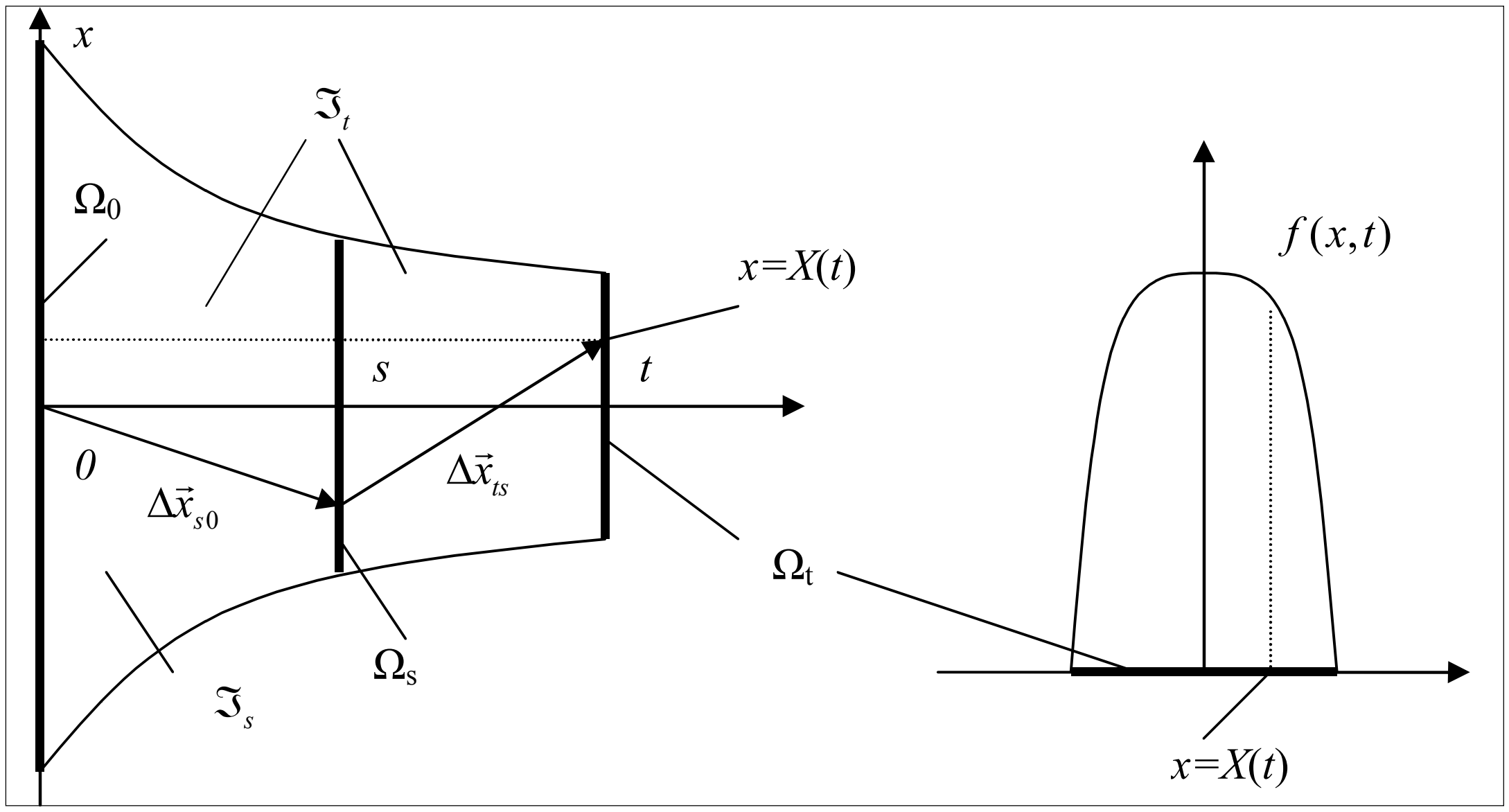

Fig. 2.3.2: Continuous $\Omega_t, \Im_t$ and $P(X_t = x \mid \Omega) = P(x \mid \Omega_t) = f(x,t)$

We can interpolate $\Omega_t$ as the snapshot of the sample space *at* time *t*, and $\Im_t$ as the set of all snapshots taken *up to* time *t*. A measurement of *X* at time *t* is taken a snapshot of sample space at time *t* as shown in Fig. 2.3.2. In this picture, we can see that $\Omega_t$ is one-dimensional front-end cross-section, while $\Im_t$ looks like a two-dimensional area.

At time *t*, the measurement picks up the value from $|\Omega_n)$, or, we can use the identity operators defined as follows:

$$I(t) = \sum_i | x_i;t)P(x_i;t| \qquad (Disacrete\, State\, SP) \qquad (2.3.5a)$$

$$I(t) = \int dx\, | x;t)P(x;t| \qquad (Continuous\, State\, SP) \qquad (2.3.5b)$$

Therefore, we can write:

$$X_t \equiv X_t(\omega) \equiv X(\omega,t) = X_t \cdot I(t) = I(t) \cdot X_t \qquad (2.3.6)$$

$$X_t \mid \Omega) = \sum_i X_t \mid x_i;t)P(x_i;t \mid \Omega) = \sum_i x_i \mid x_i)P(x_i \mid \Omega_t) \text{ (discrete } S.V.) \qquad (2.3.7a)$$

$$X_t \mid \Omega) = \int dx X_t \mid x;t)P(x;t \mid \Omega) = \int_{x\in\Omega} dx\, x \mid x)P(x \mid \Omega_t) \text{ (continuous } S.V.) \qquad (2.3.7b)$$

A *S.P. X*(*t*) has *independent increments* if, for $t_1 < t_2 < \dots < t_m < t_{m+1}$, it satisfies:

$$P(X_{t_{m+1}} - X_{t_m} = x_m \mid X_{t_{i+1}} - X_{t_i} = x_{i+1}) = P(X_{t_{m+1}} - X_{t_m} = x_m \mid \Omega),\ for\ \forall i \in \{1,\dots m-1\} \quad (2.3.8)$$

In Fig 2.3.2, we use 2-dimensional vectors to represent the displacement of *S.P X*(*t*):

$$\vec{x}_t = (t, X_t),\ \vec{x}_s = (s, X_s),\ \Delta\vec{x}_{ts} = \vec{x}_t - \vec{x}_s,\quad \Delta\vec{x}_{s0} = \vec{x}_s - \vec{x}_0 \quad (2.3.9a)$$

$$\Delta\vec{x}_{ts} \in \Im_{ts} \equiv (\Im_t - \Im_s) \cup \Omega_s,\quad \Delta\vec{x}_{s0} \in \Im_{s0} \equiv (\Im_s - \Im_0) \cup \Omega_0, \quad (2.3.9b)$$

$$\because \text{Indepedent insreaments} \Rightarrow \quad |\Im_t) = |\Im_{ts})|\Im_{s0}),$$
$$\therefore P(X_t - X_s = x \,|\, X_s - X_0 = x') = P(X_t - X_s = x \,|\, \Omega) \quad (2.3.9c)$$
$$\because \Im_t = \Im_{ts} \otimes \Im_{s0},\quad \therefore P(X_t - X_s = x \,|\, X_s - X_0 = x') = P(X_t - X_s = x \,|\, \Omega)$$

If *S.P* has independent-time increment, we can always set $X_0 = 0$ and have:

$$P(X_t = x + c \,|\, X_s = c) = P(X_t - X_s = x \,|\, X_s - X_0 = c) = P(X_t - X_s = x \,|\, \Omega) \quad (2.3.9d)$$

A *S.P* is *homogeneous* if it has the following property:

$$P(X_t - X_s = x \,|\, \Omega) = P(X_{t-\tau} - X_{s-\tau} = x \,|\, \Omega), \quad for\ t > s > \tau \geq 0 \quad (2.3.10)$$

If the *S.P.* is homogeneous and $X(0) = 0$, then we have the following property:

$$P(X_{t+s} - X_s = x \,|\, \Omega) = P(X_t - X_0 = x \,|\, \Omega) = P(X_t = x \,|\, \Omega) \equiv P(x \,|\, \Omega(t)) \quad (2.3.11a)$$

If the *S.P.* is *homogeneous with independent-time increment* and $X(0) = 0$, then we have the following property:

$$P(X_{t+s} = x + h \,|\, X_s = h) = P(X_{t+s} - X_s = x \,|\, X_s - X_0 = h) = P(X_{t+s} - X_s = x \,|\, \Omega)$$
$$= P(X_t - X_0 = x \,|\, \Omega) = P(X_t = x \,|\, \Omega) \equiv P(x \,|\, \Omega(t)) \quad (2.3.11b)$$

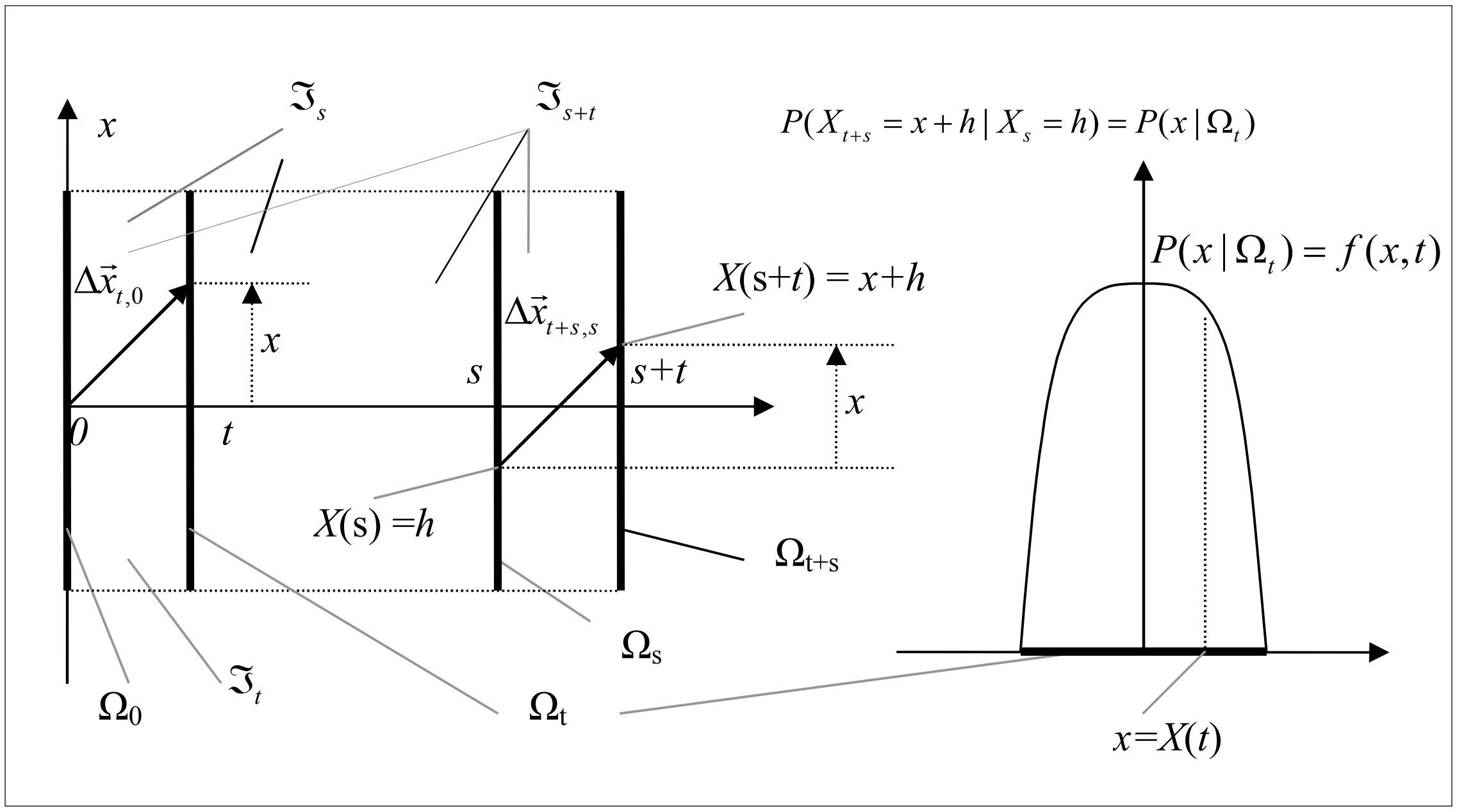

Fig. 2.3.3: Homogeneous *S.P* with Independent Increments and $X_0 = 0$

Fig 2.3.3 shows the relation between the transition probability and the absolute probability density of such a *S.P.* Because, in a Hilbert space of QM, only absolute probability density is available, any probability space induced from it should have observable(s) with property (2.3.11b).

## 2.4. Probability Space Associated with Stochastic Process

The basic notations in *PBN* for a probability space with one *S.P.* are given in [1], §5.1. Here we want to give a view which is more close to QM. To simplify our notations and reduce possible confusions, we will use the following conventions:

$$\begin{aligned}&|\Omega(t)) \equiv |\Omega^{(t)}) \equiv |\Omega_t), \quad |X(t)=x) \equiv |X_t = x) \equiv |x,t) \\ &|X_t = x, X_{t'} = x') \equiv |x,t;x',t'), \quad |X_t = x, Y_s = y) \equiv |x,t;y,s)\end{aligned} \tag{2.4.1}$$

As we have seen in Eq. (2.3.4), we can generate Eq. (2.2.6) to stochastic process to involve time evolution:

$$P:\Omega,t \mapsto |\Omega_t) \Rightarrow \text{ for } \forall x \in \Omega, \quad P_t(x) \equiv P(x,t) = P(x\,|\,\Omega_t) = f(x,t) \tag{2.4.2a}$$

This exactly is what we have for 1D QM problem in Schrodinger picture ([2], §11.12):

$$P:\Re,t \mapsto |\Psi(t)\rangle \Rightarrow \text{ for } \forall x \in \Re, \quad P(x,t) = |\langle x\,|\,\Psi(t)\rangle|^2 = |\Psi(x,t)|^2 \tag{2.4.2b}$$

In QM, we also have Heisenberg picture (see [2], §11.12) , where operators become time-dependent by using the following time-evolution unitary operator:

$$\begin{aligned}&\hat{U}(t) \equiv \exp[-it\hat{H}/\hbar], \quad \hat{U}^\dagger(t) = \hat{U}^{-1}(t), \quad |\Psi(t)\rangle = \hat{U}(t)\,|\,\Psi(0)\rangle, \quad \hat{x}(t) = \hat{U}^\dagger(t)\,\hat{x}\,\hat{U}(t), \\ &\langle\Psi(t)\,|\,\hat{x}\,|\,\Psi(t)\rangle = \langle\Psi(0)\,|\,\hat{U}^\dagger(t)\,\hat{x}\,\hat{U}(t)\,|\,\Psi(0)\rangle = \langle\Psi(0)\,|\,\hat{x}(t)\,|\,\Psi(0)\rangle\end{aligned}$$

Based on time-dependent operator, we can introduce following time-dependent states:

$$\begin{aligned}&|x(t)\rangle = \hat{U}^\dagger(t)\,|\,x\rangle, \quad \hat{U}(t)\,|\,x(t)\rangle = |\,x\rangle, \quad \langle x\,|\,\hat{U}(t) = \langle x(t)\,|, \\ &\hat{U}(t)\hat{x}(t)\,|\,x(t)\rangle = \hat{U}(t)\hat{U}^\dagger(t)\,\hat{x}\,\hat{U}(t)\hat{U}^\dagger(t)\,|\,x\rangle = \hat{x}\,|\,x\rangle = x\,|\,x\rangle\end{aligned}$$

This gives us the following important relations:

$$\langle x\,|\,\Psi(t)\rangle = \langle x\,|\,\hat{U}(t)\,|\,\Psi(0)\rangle = \langle x(t)\,|\,\Psi(0)\rangle = \Psi(x,t) \tag{2.4.2c}$$

This is common in QM. Suppose we expand the state in momentum space, we have:

$$\langle p\,|\,\Psi(t)\rangle = \langle p\,|\,\hat{U}(t)\,|\,\Psi(0)\rangle = \langle p(t)\,|\,\Psi(0)\rangle = c(p,t) \tag{2.4.2d}$$

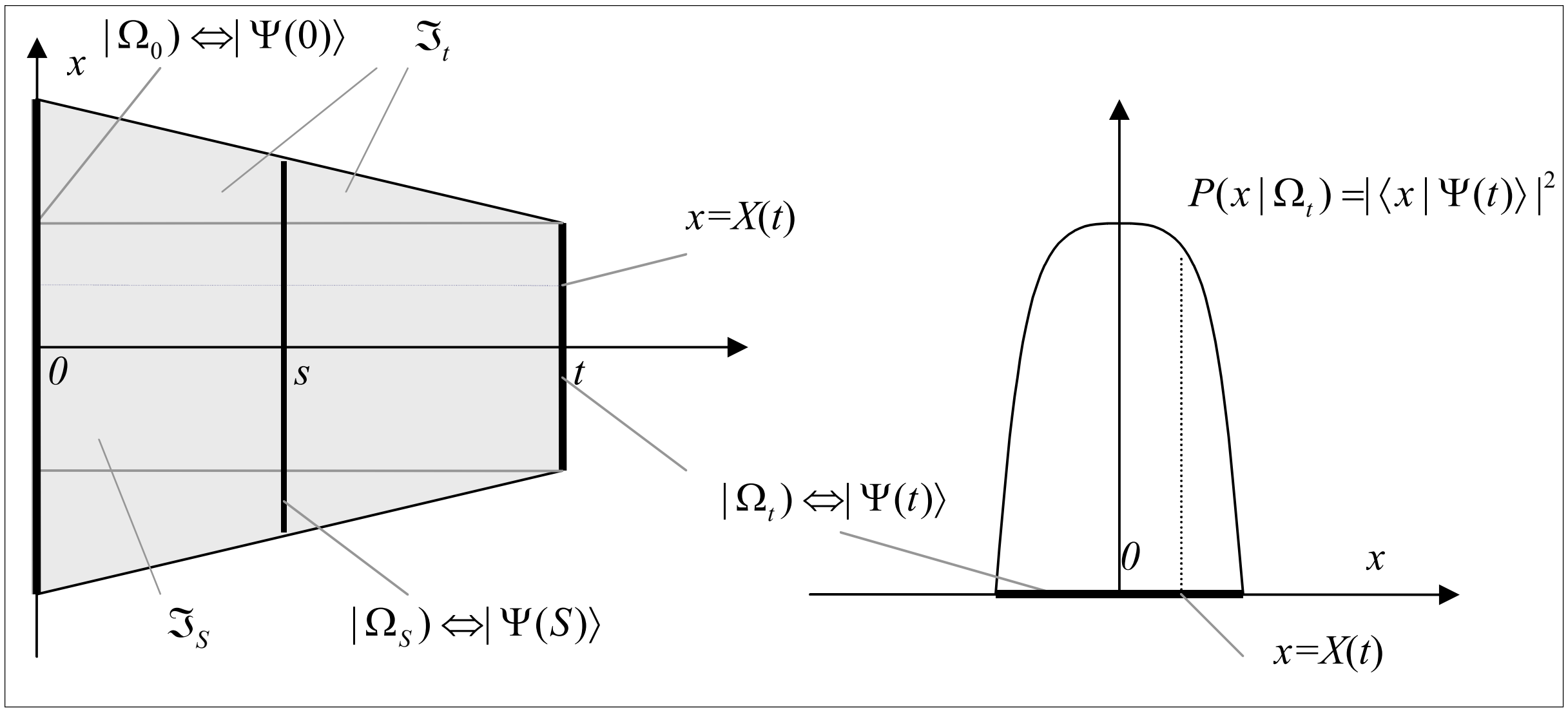

Fig.2.4.1: $|\Omega_t) \Leftrightarrow |\Psi(t)\rangle$ and $P(x|\Omega_t) \Leftrightarrow |\langle x|\Psi(t)\rangle|^2$

Fig. 2.4.1 shows the equivalence of $|\Omega_t)$ and $|\Psi(t)\rangle$ based on one *S.P X*(*t*).

Similarly, we can write the stochastic process in Heisenberg picture:

$$P:\Omega, t \mapsto |\Omega_t) \Rightarrow \text{ for } \forall x \in \Omega, \quad P_t(x) \equiv P(x,t) = P(X_t = x|\Omega) = f(x,t) \qquad (2.4.5a)$$

From Eq. (2.4.2a) and (2.4.2b), we have the following identities in *PBN*:

$$P(X_t = x|\Omega) \equiv f(x,t) \equiv P(x|\Omega_t) \qquad (2.4.5b)$$

***Dimensional Analysis***: Assume $[x]_{P.D} = L$, we have $[X_t]_{P.D} = [x]_{P.D} = L$ and:

$$[|x,t)]_{P.D} = [P(x,t|]_{P.D} = [|\Omega_0)] = L^{-1/2} \qquad (2.4.6a)$$

$$[P(\Omega|]_{P.D} = L^{1/2}, \quad [f(x,t)]_{P.D} = [P(x,t|\Omega)]_{P.D} = L^{-1} \qquad (2.4.6b)$$

If a probability space is associated with two *S.P.*, then at a given time $t \geq 0$, the base p-kets have the following properties:

$$X(t)|X_t = x, Y_t = y) \equiv X(t)|x,y;t) = x|x_i,y;t) \\ Y(t)|X_t = x, Y_t = y) \equiv Y(t)|x,y;t) = y|x_i,y;t) \qquad (2.4.7a)$$

$$P(x,y;t|x',y';t) = \delta(x-x')\delta(y-y') \quad \text{continuous state} \\ P(x,y;t|x',y';t) = \delta_{xx'}\delta_{xx'} \quad \text{discrete state} \qquad (2.4.7b)$$

$$P(x,y;t|\Omega) \equiv P(X_t = x, Y_t = y|\Omega) \equiv f(x,y;t) \equiv P(x,y|\Omega_t) \qquad (2.4.7c)$$

***Dimensional analysis***: Assuming $[x]_{P.D} = [y]_{P.D} = L$, we have:

$$[P(x,y;t|]_{P.D} = [|x,y;t)]_{P.D} = L^{-1}, \quad [f(x,y;t)]_{P.D} = L^{-2}. \tag{2.4.8a}$$

Then we can find the *P.D* for other quantities as follows:

$$|\Omega_t) = \int_{x,y} dx\,dy\,|x,y;t) f(x,y;t) \Rightarrow [|\Omega_t)]_{P.D} = L^{-1} \tag{2.4.8b}$$

$$[P(\Omega|\Omega_t)]_{P.D} = 1 \quad \Rightarrow \quad [P(\Omega|]_{P.D} = L \tag{2.4.8c}$$

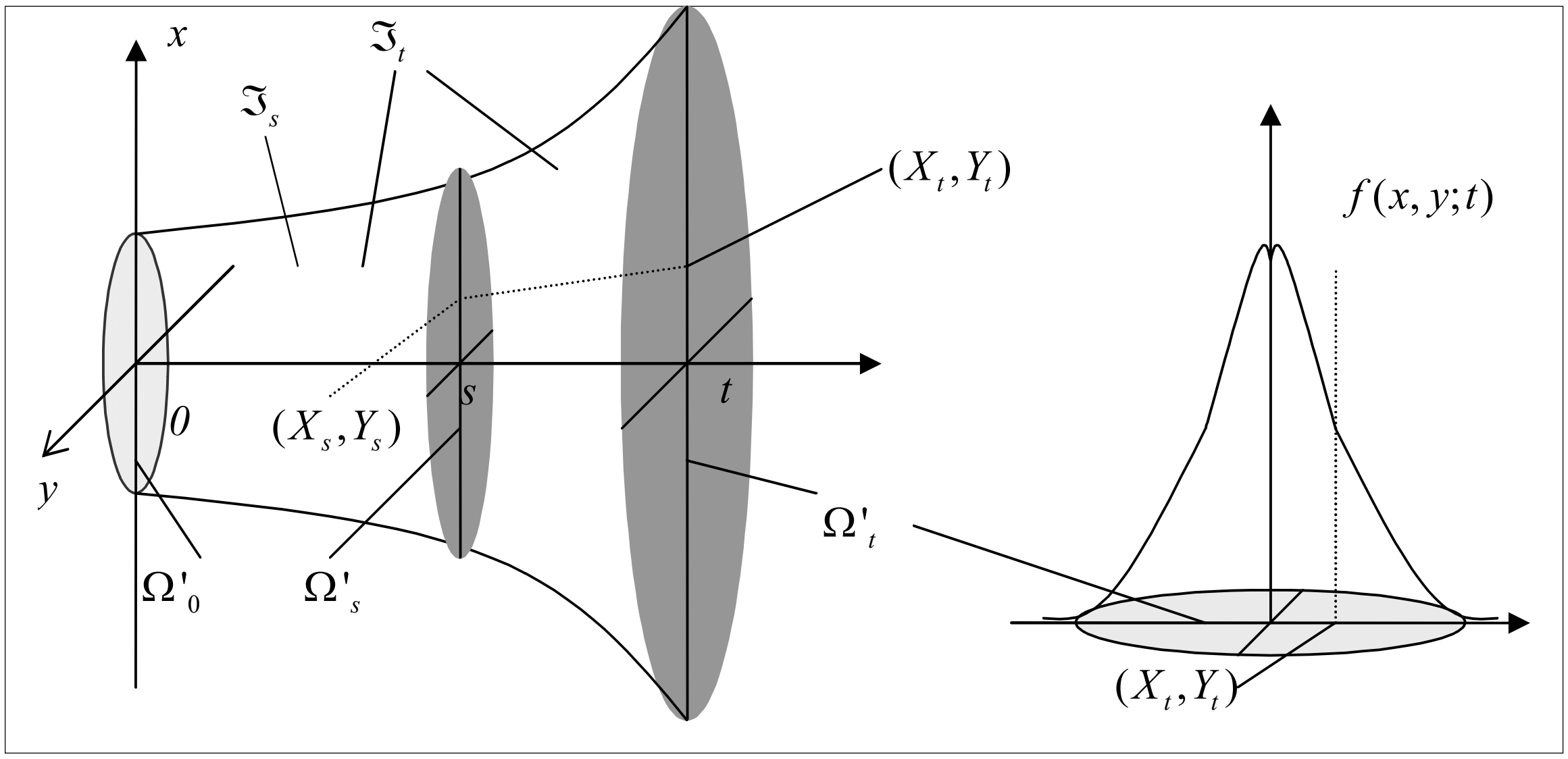

Fig.2.4.2: $\Omega'_t = \{(x,y) \in \Omega_t, P(x,y|\Omega_t) > \varepsilon\}$, $\Im_t$ and $P(x,y|\Omega_t) = f(x,y;t)$

Fig. 2.4.2 describes the relations between probability space, filtration and probability density of two *S.P*. In this picture, we define $\Omega'_t = \{(x,y) \in \Omega_t, (x,y|\Omega_t) \geq \varepsilon\}$, where $0 < \varepsilon / |\Omega_t| \ll 1$. We can see that $\Omega'_t$ is a two-dimensional front-end cross-section of $\Im_t$, which looks like a three-dimensional solid "cake".

The conditional probability of the base p-kets at different time is determined by the transition matrix element of the transition operator $\hat{P}$ in the base of v-kets. For example, for a *S.P.* $X(t)$, we may have:

$$\begin{aligned} &p(j;t+1|i;t) \equiv p_{ij}(t) = \langle i,t|\hat{P}|j;t+1\rangle \quad \text{discrete state, discrete time} \\ &p(x';t+1|x;t) \equiv \langle x;t|\hat{P}|x';t+1\rangle \quad \text{continuous state, discrete time} \end{aligned} \tag{2.4.8a}$$

$$\begin{aligned} &p(j;t|i;s) \equiv \langle i;s|\hat{P}|j;t\rangle \quad \text{discrete state, continuous time} \\ &p(x';t|x;s) = \langle x';t|\hat{P}|x;s\rangle \quad (t \geq s) \text{ continuous state, continuous time} \end{aligned} \tag{2.4.8b}$$

If both *S.P.* ***are homogeneous*** with $X(0) = 0$ and $Y(0) = 0$, then we have the following property for any $t \geq 0$ and $s \geq 0$:

$$P\left(X_{t+s}-X_s=x,Y_{t+s}-Y_s=y\,|\,\Omega\right)=P\left(X_t-X_0=x,Y_t-Y_0=y\,|\,\Omega\right)$$
$$=P\left(X_t=x,Y_t=y\,|\,\Omega\right)\equiv P(x(t),y(t)\,|\,\Omega)=P(x,y\,|\,\Omega_t)\equiv f(x,y;t) \quad (2.4.9)$$

The related transition matrices now are simplifies as:

$$p(j;t+1\,|\,i;t)=p(j;1\,|\,i;0)\equiv p(j\,|\,i)\equiv\langle i\,|\,\hat{P}\,|\,j\rangle\equiv p_{ij}\quad \text{discrete state, discrete } t$$
$$p(x';t+1\,|\,x;t)=p(x';1\,|\,x;0)\equiv\langle x;0\,|\,\hat{P}\,|\,x';1\rangle\quad \text{continuous state, discrete } t \quad (2.4.10a)$$

$$p(j;s+t\,|\,i,s)=p(j;t\,|\,i,0)\equiv\langle i;t\,|\,\hat{P}\,|\,j;0\rangle\quad \text{discrete state, continuous } t$$
$$p(x';s+t\,|\,x;s)=p(x';t\,|\,x;0)\equiv\langle x;t\,|\,\hat{P}\,|\,x;0\rangle\quad \text{continuous state, continuous } t \quad (2.4.10b)$$

A *S.P X*(t) may have ***Markov property***, which assumes that the future probability distribution can be predicted from the current system state, but not the past system state. This means, for $t_1 < t_2 < t_{\dots} < t_m < t_{m+1}$,

$$P(X(t_{m+1})=x_{m+1}\,|\,X(t_m)=x_m,X(t_{m-1})=x_{m-1},\dots,X(t_1)=x_1)$$
$$\equiv P(x_{m+1},t_{m+1}\,|\,x_m,t_m;x_{m-1},t_{m-1};\dots;x_1,t_1)=P(x_{m+1},t_{m+1}\,|\,x_m,t_m) \quad (2.4.11)$$

We can generate Eq. (2.3.5b) to two *S.P*. of continuous-states:

$$I(t)=\int_{x,y\in\Omega}dx\,dy\,|\,x,y;t)\,p(x,y;t\,| \quad (2.4.12)$$

It is time dependent, because when we insert into a p-bracket (a conditional probability by definition), the newly formed p-bracket may be a transition matrix if its p-ket and p-bra have different time), and it may vanish if time is not properly chosen. One good example is the *Chapman-Kolmogorov theorem*. This equation can be derived by using Conditional Total Probability Law (TPL) and Markov property. Now, if a *S.P.* has *Markov property*, we can "derive" them in *PBN* by simply inserting our identity operator Eq. (2.4.7b) with appropriate time:

$$p^{m+n}{}_{ij}\equiv p(j;m+n\,|\,i;0)=\sum_k p(j;m+n\,|\,k;m)p(k;m\,|\,i;0)$$

$$=\ \sum_k p^m{}_{ik}\ p^n{}_{kj}\qquad (\textit{Disacrete time, discrete R.V}) \quad (2.4.13)$$

$$p_{ij}(t+s)\equiv p(j;t+s\,|\,i;0)=\sum_k p(j;t+s\,|\,k;s)p(k;s\,|\,0;i)$$

$$=\ \sum_k p_{ik}(t)\ p_{kj}(s)\qquad (\textit{Continuous time, discrete R.V}) \quad (2.4.14)$$

$$p(x;m+n\,|\,x';0)=\int dx''p(x;m+n\,|\,x'';m)p(x'';m\,|\,0;x')$$
$$(\textit{diccrete time, continuous R.V}) \quad (2.4.15)$$

$$p(x',y';t\,|\,x,y,s) = p(x',y';t\,|\,I(\tau)\,|\,x,y,s) = \int dx''dy''\,p(x',y';t\,|\,x'',y'';\tau)\,p(x'',y'';\tau\,|\,x,y;s)$$
$$\text{where } t > \tau > s \qquad (Continuous\,time, two\,continuous\,R.V) \qquad (2.4.16)$$

In general, if a *S.P.* has Markov property, then we can insert an Identity operator inside the transition matrix (a p-bracket), with a time less than the time on the left and greater than the time on the right.

*Eigenvectors of Homogeneous Markov Chain*: Assuming we have a homogeneous MC (Markov chain) and its Transition operator $\hat{P}$ has a right eigenvector $|\varphi\rangle$ (a v-ket) corresponding to eigenvalue $\lambda$, then we can simply write:

$$\hat{P}\,|\,\varphi\rangle = \lambda\,|\,\varphi\rangle \qquad (2.4.17)$$

If the MC has discrete time increment and discrete states, inserting the unit operator in vector space, we can write:

$$\langle i\,|\,\hat{P}\,|\,\varphi\rangle = \sum_j \langle i\,|\,\hat{P}\,|\,j\rangle\langle j\,|\,\varphi\rangle = \sum_j p_{ij}\langle j\,|\,\varphi\rangle = \sum_j p_{ij}\varphi(j) = \langle i\,|\,\lambda\,|\,\varphi\rangle \qquad (2.4.18a)$$

$$\text{Or: } \sum_k p_{ij}\varphi(j) = \lambda\varphi(i) \qquad (2.4.18b)$$

An element of transition matrix of discrete-state *M.C* is a dimensionless probability by nature. From Eq. (2.4.18b), we have following dimensional analysis:

$$\text{Discret-state } M.C\text{: } [P]_{P.D} = [P_{ij}]_{P.D} = [\lambda]_{P.D} = 1 = L^0 \qquad (2.4.18c)$$

If the MC has discrete time increment and continuous states, then the transition matrix is expressed as $\langle y\,|\,\hat{P}\,|\,z\rangle \equiv p(z,1\,|\,y,0)$. Inserting the unit operator in (2.2.7a), we have:

$$\lambda\langle y\,|\,\varphi\rangle = \langle y\,|\,\hat{P}\,|\,\varphi\rangle = \langle y\,|\,\hat{P}\cdot I\,|\,\varphi\rangle = \langle y\,|\,\hat{P}\int dz\,|\,z\rangle dz\,\langle z\,|\,\varphi\rangle \qquad (2.4.19a)$$

$$\text{Or: } \lambda\varphi(y) = \int dz\langle y\,|\,\hat{P}\,|\,z\rangle\langle z\,|\,\varphi\rangle = \int dz\,p(z,1\,|\,y,0)\varphi(z) \qquad (2.4.19b)$$

Assuming $[x]_{P.D} = L$, and using the physical dimensions given in Eq. (2.4.6), we have following dimensional analysis from Eq. (2.4.19b):

$$\text{Continuous-state } M.C\text{: } [\langle x\,|\,\hat{P}\,|\,z\rangle]_{P.D} = L^{-1}, \quad [\hat{P}]_{P.D} = [\lambda]_{P.D} = 1 \qquad (2.4.20)$$

## 3. Properties of Conditional Expectation in *PBN*

In this section, we investigate properties of Conditional Expectation (*C.E*) using *PBN*, as the preparation of introducing martingales.

## 3.1. The Properties of Conditional Expectation Related to *R.V.*

We discuss here a few important properties of conditional expectation (*CE*) related to *R.V*. Some of them are going to be used in investigating martingales based on sequence of *R.V.*

**Conditional Expectation (*CE*)**: Let $X$ be a $R.V$ on $(\Omega, \Im, P)$, $H \in \Im$ and $\mathrm{P}(H \mid \Omega) > 0$; then the conditional probability density can be repressed as follows for *discrete states*:

$$E(X \mid H) \equiv P(\Omega \mid X \mid H) = \sum_i P(\Omega \mid X \mid x_i) P(x_i \mid H) = \sum_i x_i P(x_i \mid H) \tag{3.1.1a}$$

$$P(x_i \mid H) = \frac{P(x_i \cap H \mid \Omega)}{P(H \mid \Omega)} = \begin{cases} P(x_i \mid \Omega) & if\ x_i\ independent\ of\ H \\ \boldsymbol{1}_H(x_i) \dfrac{P(x_i \mid \Omega)}{P(H \mid \Omega)} & otherwise \end{cases} \tag{3.1.1b}$$

$$P(\Omega \mid X \mid H) = \sum_i x_i P(x_i \mid H) = \frac{\sum_{x_i \in H} x_i P(x_i \mid \Omega)}{P(H \mid \Omega)} \tag{3.1.1c}$$

For continuous states and (see Eq. (1.1.1 1b) and (1.1.1 4b)) we have similar definitions:

$$E(X \mid H) \equiv P(\Omega \mid X \mid H) = \int P(\Omega \mid X \mid x) dx\, P(x \mid H) = \int x\, dx\, P(x \mid H) \tag{3.1.2a}$$

where $P(x \mid H) = \dfrac{P(x \cap H \mid \Omega)}{P(H \mid \Omega)}$ if $[P(H \mid \Omega)]_{P.D} = 1$ and $x \in H$ (3.1.2b)

Otherwise, if $x$ is independent of $H : P(x \mid H) = P(x \mid \Omega)$ (3.1.2c)

if $H = x' \in \Omega : P(x \mid x') = \delta(x - x')$ (3.1.2d)

We can derive some interesting expressions (see Eq. (3.2.5-3.2.7) in [4]) related to expectation value using indicator function defined in (2.1.10):

1. $\langle \boldsymbol{1}_B \rangle \equiv P(\Omega \mid \boldsymbol{1}_B \mid \Omega) \underset{Eq.(2.1.12)}{=} P(B \mid \Omega)$. (3.1.3a)
2. $P(\Omega \mid X\boldsymbol{1}_B \mid \Omega) = P(B \mid \Omega) P(\Omega \mid X \mid B),\ where\ P(B \mid \Omega) > 0$. (3.1.3b)

   *Proof*: It is trivial for discrete states. But for *continuous case*, as mentioned in §3.2 of Ref. [4], that the proof needs to use Measure theory. Our following proof in *PBN* seems not to need that:

$$P(\Omega \mid X\boldsymbol{1}_B \mid \Omega) \underset{(2.1.10)}{=} \int_{x \in B} dx\, P(\Omega \mid X \mid x) P(x \mid \Omega)$$

$$= \int_{x \in B} dx\, P(\Omega \mid x \mid x) P(x \mid \Omega) = \int_{x \in B} dx\, x\, P(x \mid \Omega) = P(B \mid \Omega) \frac{\int_{x \in \Omega} dx\, x\, P(x \mid \Omega)}{P(B \mid \Omega)}$$

$$\underset{(3.1.2b)}{=} P(B \mid \Omega) P(\Omega \mid X \mid B) \text{ (for continuous states)}$$

3. $P(\Omega \mid \boldsymbol{1}_A \mid B) = P(A \mid B), \textit{ where } P(B \mid \Omega) > 0$. (3.1.4a)

$$proof: P(B \mid \Omega) P(\Omega \mid \boldsymbol{1}_A \mid B) \underset{(3.1.3b)}{=} P(\Omega \mid \boldsymbol{1}_A \boldsymbol{1}_B \mid \Omega) = P(\Omega \mid \boldsymbol{1}_{A \cap B} \mid \Omega)$$
$$\underset{(3.1.3a)}{=} P(A \cap B \mid \Omega) \underset{(3.1.2b)}{=} P(B \mid \Omega) P(A \mid B) \quad (3.1.4b)$$

It is easy to see that, for a Borel function $g(x)$, we have for discrete states:

$$E(g(X) \mid H) \equiv P(\Omega \mid g(X) \mid H) = \sum_i P(\Omega \mid g(X) \mid x_i) P(x_i \mid H) = \sum_i g(x_i) P(x_i \mid H) \quad (3.1.5a)$$

And for continuous states:

$$E(g(X) \mid H) \equiv P(\Omega \mid g(X) \mid H) = \int dx\, P(\Omega \mid g(X) \mid x) P(x \mid H) = \int dx\, g(x) P(x \mid H) \quad (3.1.5b)$$

If $X$ and $Y$ are *two observables* on a probability space $(\Omega, \Im, P)$, and $g(x)$ a Borel function, then the *C.E* of $g(X)$ given $Y = y$ can be written by using Eq. (2.2.8b) (see also [4], §3.2):

$$E[g(X) \mid Y = y] \equiv P(\Omega \mid g(X) \mid y) = \int P(\Omega \mid g(X) \mid x)\, dx\, P(x \mid y)$$
$$= \int g(x) dx P(x \mid y) \quad (3.1.6a)$$

$$where: P(x \mid y) = P(x, * \mid *, y) = \frac{P(X = x, Y = y \mid \Omega)}{P(*, Y = y \mid \Omega)} = \frac{P(x, y \mid \Omega)}{P(y \mid \Omega)} \equiv \frac{f(x, y)}{P(y \mid \Omega)} \quad (3.1.6b)$$

If $X$ and $Y$ are independent, by using Eq. (2.2.8c) and (2.2.12), we have:

$$P(x \mid y) = P(x, * \mid *, y) = \frac{P(x, y \mid \Omega)}{P(*, y \mid \Omega)} = \frac{P(x \otimes y \mid \Omega)}{P(y \mid \Omega_Y)}$$
$$= \frac{P(x \mid \Omega_X) P(y \mid \Omega_Y)}{P(y \mid \Omega_Y)} = P(x \mid \Omega_X) \quad (3.1.6c)$$

Dimensional Analysis shows our expression is consistent (see also Eq. (2.2.13)):

$$[P(x, y \mid y)]_{P.D} = [P(x \mid y)]_{P.D} = [P(x, * \mid \Omega)]_{P.D} = [P(x \mid \Omega_X)]_{P.D} = L^{-1}.$$

Note: please do not confuse the conditional probability $P(x \mid y)$ of two *R.V.* ($X$ and $Y$) with the conditional probability of one *R.V.* ($X$), where $P(x \mid x') \equiv P(X = x \mid X = x') = \delta(x - x')$.

Now we list some important *C.E* properties with proofs in *PBN* for continuous cases.

***CE Property 1***: If $g \geq 0$, then from Eq. (3.1.6) and the fact that $P(x \mid y) \geq 0$, we have $E[g(X) \mid y] \geq 0$ for all $y \in \Im$. Replace y with Y, we have:

If $g \geq 0$, then $E[g(X) \mid Y] \geq 0$ (3.1.7)

***CE Property 2***:
$$P(\Omega \mid [a_1 g_1(X_1) + a_2 g_2(X_2)] \mid Y) = a_1 P(\Omega \mid g_1(X_1) \mid Y) + a_2 P(\Omega \mid g_2(X_2) \mid Y) \quad (3.1.8)$$

*Proof*: For $\forall y \in \Im_Y$,
$$P(\Omega \mid [a_1 g_1(X_1) + a_2 g_2(X_2)] \mid y) = \int dx_1 dx_2 P(\Omega \mid [a_1 g_1(X_1) + a_2 g_2(X_2)] \mid x_1, x_2) P(x, x_2 \mid y)$$
$$= \int dx_1 a_1 g_1(x_1) \int dx_2 P(x_1, x_2 \mid y) + \int dx_1 dx_2 a_2 g_2(x_2) P(x_1, x_2 \mid y)$$

$$= a_1 \int dx_1 g_1(x_1) P(x_1 \mid y) + a_2 \int dx_2 g_2(x_2) P(x_2 \mid y) = a_1 P(\Omega \mid g_1(X_1) \mid y) + a_2 P(\Omega \mid g_2(X_2) \mid y)$$

***CE Property 3***: If *X* and *Y* are independent, we have for any $y \in \Im$:

$$P(\Omega \mid g(X) \mid y) = \int P(\Omega \mid g(X) \mid x)\, dx\, P(x \mid y) = \int g(x) dx\ P(x \mid \Omega) = \langle g(X) \rangle$$
$$\therefore \quad P(\Omega \mid g(X) \mid Y) = \langle g(X) \rangle \quad (3.1.9)$$

***CE Property 4***: $\langle (\Omega \mid g(X) \mid Y) h(Y) \rangle = \langle g(X) h(Y) \rangle$ (3.1.10)

Proof: $\langle (\Omega \mid g(X) \mid Y) h(Y) \rangle = \int h(y) P(\Omega \mid g(X) \mid y)\, dy\, P(y \mid \Omega_Y) =$
$$\iint h(y)\, dx\, g(x) P(x \mid y) dy\, P(y \mid \Omega) \underset{(3.1.6b)}{=} \iint h(y)\, dx\, g(x) \frac{P(x, y \mid \Omega)}{P(y \mid \Omega)} dy\, P(y \mid \Omega)$$
$$= \iint dx\, dy\, g(x)\, h(y)\, P(x, y \mid \Omega) = \langle g(X)\, h(Y) \rangle \quad (3.1.11)$$

***CE Property 5*** (*Double Expectation Formula)*: $\langle P(\Omega \mid g(X) \mid Y) \rangle = \langle g(X) \rangle$ (3.1.12a)
Proof: This is a special case of property 4, when we have $h(y) = 1\ for\ \forall y \in \Im$.
This can be extended to: $\langle P(\Omega \mid g(X) \mid Y_1, \cdots Y_n) \rangle = \langle g(X) \rangle$ (3.1.12b)

***CE Property 6***: $P(\Omega \mid g(X) h(Y) \mid Y) = h(Y) P(\Omega \mid g(X) \mid Y)$ (3.1.13)
Proof: For $\forall y \in \Im$: $P(\Omega \mid g(X) h(Y) \mid y) = P(\Omega \mid g(X) h(y) \mid y) = h(y) P(\Omega \mid g(X) \mid y)$
Replace y with Y and we have the proof.

***CE Property*** *7*: Set $g(X) = 1$ in Eq. (3.1.13), we have: $P(\Omega \mid h(Y) \mid Y) = h(Y)$ (3.1.14a)
Similarly, we have for *n R.V*'s: $P(\Omega \mid h(Y_1, ..., Y_n) \mid Y_1, ..., Y_n) = h(Y_1, ..., Y_n)$ (3.1.14b)

***CE Property 8***: $P(\Omega \mid [P(\Omega \mid g(X) \mid Z)] \mid Y, Z) = P(\Omega \mid g(X) \mid Z)$ (3.1.15)
Proof: Defining $h(Y, Z) = P(\Omega \mid g(X) \mid Z)$, then use Eq. (3.1.14b) with $Y_1 \to X, Y_2 \to Z$.

***CE Property 9***: $P(\Omega \mid [P(\Omega \mid g(X) \mid Y, Z)] \mid Z) = P(\Omega \mid g(X) \mid Z)$ (3.1.16)

Proof: Defining $h(Y,Z) = P(\Omega\,|\,g(X)\,|\,Y,Z)$, then we have:

$$P(\Omega\,|\,h(Y,Z)\,|\,z) = P(\Omega\,|\,h(Y,z)\,|\,z) = \int h(y,z)\,dy\,P(y\,|\,z)$$

$$= \int P(\Omega\,|\,g(X)\,|\,y,z)\,dy\,P(y\,|\,z) = \iint g(x)\,dx\,P(x\,|\,y,z)\,dy\,P(y\,|\,z) \qquad (3.1.17a)$$

But: $\int P(x\,|\,y,z)\,dy\,P(y\,|\,z) \underset{(3.1.6b)}{=} \int dy \dfrac{P(x,y,z\,|\,\Omega)}{P(y,z\,|\,\Omega)} \dfrac{P(y,z\,|\,\Omega)}{P(z\,|\,\Omega)}$

$$= \frac{\int dy P(x,y,z\,|\,\Omega)}{P(z\,|\,\Omega)} \equiv \frac{P(x,z\,|\,\Omega)}{P(z\,|\,\Omega)} \underset{(3.1.6b)}{=} P(x\,|\,z) \qquad (3.1.17b)$$

Inserting (3.1.17b) into (3.1.17a), we have:

$$P(\Omega\,|\,h(Y,Z)\,|\,z) = \int g(x)\,dx\,P(x\,|\,z) = P(\Omega\,|g(X)\,|\,z) \qquad \text{QED}$$

Replacing $Z \Rightarrow Y_1,\cdots Y_n$ and $Y \Rightarrow Y_{n+1}$, Eq. (3.1.16) can be rewritten to:

$$P(\Omega\,|\,[P(\Omega\,|\,g(X)\,|\,Y_1,\cdots Y_n,Y_{n+1})]\,|\,Y_1,\cdots Y_n) = P(\Omega\,|\,g(X)\,|\,Y_1,\cdots Y_n) \qquad (3.1.18)$$

## 3.2. The Properties of Conditional Expectation Based on Sigma-fields

In this section, we will represent the properties of *CE* based on σ-fields, because some of them will be used for martingales based on filtration latter. We will only provide some hints or shortcut of their proof. First, we need to introduce some definitions related to *R.V.* and σ-fields.

**Definition 3.2.1** ([4], definition 9.5.1; [5], definition 1.2.8): Given a *R.V X* we denote by σ (*X*) the smallest σ-field $\mathfrak{B} \subseteq \mathfrak{I}$ such that *X*(*ω*) is measurable on (Ω, $\mathfrak{B}$ ). We call the σ-field generated by *X* and denote it σ (*X*) or $\mathfrak{I}_X \equiv \mathfrak{I}(X)$.

**Definition 3.2.2** ([4], definition 9.5.6): Assume that *X* is *R.V* on probability space $(\Omega,\mathfrak{I},P)$, $\langle|\,X\,|\rangle < +\infty$, and σ-field $\mathrm{B} \subset \mathfrak{I}$, then the *conditional expectation* of *X* about $\mathfrak{B}$ is a *R.V.*: $Z(X) = (\Omega\,|\,X\,|\,\mathrm{B})$ satisfying following conditions:

(1). *Z*(*X*) is measurable on $\mathfrak{B}$: $\{\omega, Z(X) \le x\} = \{\omega, (\Omega\,|\,X\,|\,\mathrm{B}) \le x\} \in \mathrm{B}$

(2). For $\forall B \in \mathrm{B}$, $\langle X\boldsymbol{1}_B\rangle = P(\Omega\,|\,X\boldsymbol{1}_B\,|\,\Omega) = \langle P(\Omega\,|\,X\,|\,\mathrm{B})\boldsymbol{1}_B\rangle$.

*Note*: From Eq. (3.1.3b), $P(\Omega\,|\,X\boldsymbol{1}_B\,|\,\Omega) = P(B\,|\,\Omega)P(\Omega\,|\,X\,|\,B),\ where\ P(B\,|\,\Omega) > 0$

According to [4], by using measure theory, the second condition is equivalent to:
(2'). For any *R.V Y*, bounded and measurable on $\mathfrak{B}$: $P(\Omega\,|\,XY\,|\,\Omega) = \langle P(\Omega\,|\,X\,|\,\mathrm{B})Y\rangle$.

**Definition 3.2.3** ([4], definition 9.5.7; [5]): If *X* and *Y* are *R.V* on probability space $(\Omega,\mathfrak{I},P)$ and $\langle|\,X\,|\rangle < +\infty$, then the *C.E* of *X* given *Y* is defined as:

$$P(\Omega \mid X \mid Y) \equiv P(\Omega \mid X \mid \Im(Y)) \qquad (3.2.1a)$$

Here $\Im(Y)$ is the σ-field generated by *Y*. The following definition is equivalent to Definition 3.2.3:

**Definition 3.2.4** ([5], definition 2.1.5): The *conditional expectation* of $X \in L^1(\Omega, \Im, P)$ given a σ-field $\mathsf{B} \subseteq \Im$, is the *R.V Z* on $(\Omega, \mathsf{B})$ such that

$$P(\Omega \mid (X - Z)I_B \mid \Omega) = 0, \quad \forall B \in \mathsf{B} \subseteq \Im \qquad (3.2.1b)$$

Therefore, $P(\Omega \mid X \mid Y)$ corresponds to the special case of $\mathsf{B} = \Im(Y) = \Im_Y$ as in (3.2.1a).

**Example 3.2.1** ([5], Example 2.3.1): If $\mathcal{B}$ and $\sigma(X)$ are independent (i.e. for all $A \in \mathsf{B}$ and $B \in \sigma(X)$ we have that $P(\Omega \mid A \cap B) = P(\Omega \mid A)P(\Omega \mid B)$, then, $\langle XI_A \rangle = \langle X \rangle \langle I_A \rangle$ for all $A \in \mathsf{B}$ . So, the constant $Z = P(\Omega \mid X \mid \Omega) = \langle X \rangle$ satisfies (3.2.1b), that is, in this case $P(\Omega \mid X \mid \mathsf{B}) = P(\Omega \mid X \mid \Omega) = \langle X \rangle$ .

**Example 3.2.2** ([5], Example 2.3.2): Suppose $X \in L^1(\Omega, \Im, P)$ . Obviously, $Z = X$ satisfies (3.2.1b) and by our assumption *X* is measurable on $\mathcal{B}$. Consequently, here we see that $P(\Omega \mid X \mid \mathsf{B}) = X$ almost sure. In particular, if *X* = c almost sure (a constant *R.V.*), then we have that $P(\Omega \mid X \mid \mathsf{B}) = c$ almost sure.

We see that if *X is bounded and measurable on* $\mathcal{B}$, then for $\forall x \in \mathsf{B} \subseteq \Im$ we have:

$$P(\Omega \mid X \mid x \in \mathsf{B}) = P(\Omega \mid x \mid x) = x\,P(\Omega \mid x) \underset{(1.1.8b)}{=} x \iff P(\Omega \mid X \mid \mathsf{B}) = X \qquad (3.2.1c)$$

The *C.E* defined in Eq. (3.2.1a or b) has the following properties (see [4], §9.5.2 and [5], §2.3). They can be easily derived by directly using Eq. (3.2.1c) and *PBN* as in hints, or by replacing the condition from a *R.V* to $\mathcal{B}$ in properties of *CE* in §3.1, as in shortcuts.

(1). [Similar to ***CE Property 2***] if $X_1$ and $X_2$ are bounded and measurable on $\mathcal{B}$

$$P(\Omega \mid (a_1X_1 + a_2X_2) \mid \mathsf{B}) = a_1P(\Omega \mid X_1 \mid \mathsf{B}) + a_2P(\Omega \mid X_2 \mid \mathsf{B}) \qquad (3.2.2)$$

Hint: as linear combination of operators, $(a_1X_1 + a_2X_2) \mid \mathsf{B}) = a_1X_1 \mid \mathsf{B}) + a_2X_2 \mid \mathsf{B})$

Shortcut: replace $|\,Y)$ with $|\,\mathsf{B}\,)$ in Eq. (3.1.8).

(2). [Similar to ***CE Property 1***] if *X* is bounded and measurable on $\mathcal{B}$

If $X \ge 0$, then $\mathrm{P}(\Omega \mid X \mid \mathsf{B}) \ge 0$ (3.2.3)

Hint: $P(\Omega \mid X \mid x \in \mathsf{B}) = P(\Omega \mid x \mid \mathsf{B}) = x > 0$

Shortcut: replace $|\,Y)$ with $|\,\mathsf{B}\,)$ in Eq. (3.1.7).

(3). [Similar to ***CE Property 6***] for every *R.V. Y*, bounded and measurable on $\mathcal{B}$,

$$P(\Omega \mid XY \mid \mathsf{B}) = Y\,P(\Omega \mid X \mid \mathsf{B}) \qquad (3.2.4)$$

Hint: $P(\Omega \mid XY \mid y \in \mathrm{B}) = P(\Omega \mid Xy \mid \mathrm{B}) = y\,P(\Omega \mid X \mid \mathrm{B})$
Shortcut: replace $|Y)$ with $|\mathrm{B}\ )$ in Eq. (3.1.13).

(4). [Similar to ***CE Property 4***] for every *R.V. Y*, bounded and measurable on $\mathfrak{B}$,

$$\langle XY\rangle = \langle Y\,P(\Omega \mid X \mid \mathrm{B})\rangle \qquad (3.2.5)$$

*Proof*: $\langle P(\Omega \mid X \mid \mathrm{B})Y\rangle = \int dy P(\Omega \mid X \mid y \in \mathrm{B})y\,P(y \mid \Omega) = \iint dx\,dy\,xy\,P(x \mid y)P(y \mid \Omega)$

$$\underset{(3.1.6b)}{=} \iint dx\,dy\,x\,yP(x, y \mid \Omega) = \langle XY\rangle \qquad \text{QED}$$

Shortcut: replace $|Y)$ with $|\mathrm{B}\ )$ in Eq. (3.1.10). Note: Eq. (3.2.5) is the condition (2') in ***Definition 3.2.2***, and we derive it without directly using Measure theory.

(5). [Similar to ***CE Property 7***] for every *R.V. Z*, bounded and measurable on $\mathfrak{B}$,

$$P(\Omega \mid Z \mid \mathrm{B}) = Z \qquad (3.2.6)$$

Hint: same as Eq. (3.2.1).
Shortcut: replace $|Y)$ with $|\mathrm{B}\ )$ and let $Z = h(Y)$ in Eq. (3.1.14a).

(6). ***Tower Property*** [see [7], lemma 5.2): If σ-field $\mathcal{A} \subseteq \mathrm{B} \subseteq \Im$, then

$$P(\Omega \mid [P(\Omega \mid X \mid \mathrm{B})] \mid \mathcal{A}) = P(\Omega \mid X \mid \mathcal{A}) \qquad (3.2.7a)$$

$$P(\Omega \mid [P(\Omega \mid X \mid \mathcal{A})] \mid \mathrm{B}) = P(\Omega \mid X \mid \mathcal{A}) \qquad (3.2.7b)$$

Hint for (3.2.7a): $P(\Omega \mid [P(\Omega \mid X \mid x \in \mathcal{A})] \mid \mathcal{A}) = P(\Omega \mid x \mid \mathcal{A})$
Hint for (3.2.7b): $P(\Omega \mid [P(\Omega \mid X \mid x \in \mathcal{A})] \mid \mathrm{B}) = P(\Omega \mid (x \in \mathcal{A}) \mid \mathrm{B}) = (\Omega \mid x \mid \mathcal{A})$

(7). [Similar to ***CE Property 5***:] for every *R.V. X*, bounded and measurable on $\mathfrak{B}$,

$$\langle X\rangle = \langle P(\Omega \mid X \mid \mathrm{B})\rangle \qquad (3.2.8)$$

Hint: take expectation on both sides of (3.2.6)
Shortcut: replace $|Y)$ with $|\mathrm{B}\ )$ in Eq. (3.1.12).

(8). If a *R.V. X* is independent on $\mathfrak{B}$, then by using Eq. (3.1.9), we have

$$P(\Omega \mid X \mid \mathrm{B}) = P(\Omega \mid X \mid \Omega) = \langle X\rangle \qquad (3.2.9)$$

Hint: this is the same as **Example 3.2.1**.

In conclusion, if a property is true for a *C.E* given |*Y*), then it is true for the same *C.E* given $|\mathrm{B}) = |\Im_Y)$. Or, we have the following equivalence in *C.E* properties: For *R.V. X*, adapted to filtration $|\mathrm{B}) = |\Im_X)$, we have for a Borel function *g*(*x*):

$$P(\Omega \mid g(X) \mid \mathrm{B}) = P(\Omega \mid g(X) \mid \Im_X) = g(X) \;\Leftrightarrow\; P(\Omega \mid g(X) \mid X) = g(X) \qquad (3.2.10)$$

We can easily extend this property to *S.P*. Based on discussion in §2.3 and §3.1, we have

$$P(\Omega \mid g(X_t) \mid \mathrm{B} = \Im_s) = P(\Omega \mid g(X_t) \mid \Im_s) = g(X_s) \text{ for all } t \le s \qquad (3.2.11)$$

$$P(\Omega \mid g(X_t) \mid \mathrm{B} = \Im_s) = P(\Omega \mid g(X_t) \mid \Im_s) = \langle g(X_t)\rangle \text{ for all } t > s \qquad (3.2.12)$$

For a *S.P.* with *independent time increments*, from Eq. (2.3.8), we have, for any $t > s$:

$$P(\Omega \mid \{X_t - X_s\} \mid \Im_s) = P(\Omega \mid \{X_t - X_s\} \mid \Omega) = \langle X_t - X_s \rangle \quad (3.2.13)$$

$$\text{If also homogeneous}: \ \mathrm{P}(\Omega \mid \{X_t - X_s\} \mid \Im_s) = \langle X_{t-s} \rangle = \int_{x \in \Omega} x\, dx f(x; t-s) \quad (3.2.14)$$

With the help of Eq. (3.2.10-14), we can readily apply the above *C.E* properties *C.E* given filtration as condition, which will be used in our next topic on martingales. Again, our proofs in *PBN* seem not to need direct usage of measure theory.

## 4. Introduction to Martingales

In this section, we first introduce the definitions of martingale based on random sequence (more likely, from stochastic processes of discrete time) and some examples. Then we introduce martingales define on σ-field. The contents of this section mainly come from Chapter 9, Ref. [4], Chapter 4 of Ref [5] and Lecture 3 of Ref [6].

### 4.1. Discrete-Time Martingale based on *R.V.* Sequences

Here are the definitions of martingales (denoted *M.G.*) given by Ref. [8].
**Definition 4.1.1**: A discrete-time **martingale** is a discrete-time stochastic process (i.e., a sequence of random variables) $X_1, X_2, X_3 \ldots$ that satisfies for all $n$

$$\begin{aligned} &E(|X_n|) \equiv \langle |X_n| \rangle = P(\Omega \mid \{|X_n|\} \mid \Omega) < \infty \\ &E(X_{n+1} \mid X_1, \ldots, X_n) \equiv P(\Omega \mid X_{n+1} \mid X_1, \ldots, X_n) = X_n \end{aligned} \quad (4.1.1)$$

i.e., the conditional expected value of the next observation, given all of the past observations, is equal to the last observation.

**Definition 4.1.2**: Somewhat more generally, a sequence $Y_1, Y_2, Y_3 \ldots$ is said to be a martingale with respect to another sequence $X_1, X_2, X_3 \ldots$ if for all $n$:

$$E(|Y_n|) \equiv \langle |Y_n| \rangle \equiv P(\Omega \mid \{|Y_n|\} \mid \Omega) < \infty \quad (4.1.2a)$$

$$E(Y_{n+1} \mid X_1, \ldots, X_n) \equiv P(\Omega \mid Y_{n+1} \mid X_1, \ldots, X_n) = Y_n \quad (4.1.2b)$$

From Eq. (4.1.2), we see that $Y_n$ is a function of $X_1, X_2, X_3 \ldots$, that is:

$$Y_n = P(\Omega \mid Y_{n+1} \mid X_1, \ldots, X_n) \quad (4.1.3)$$

Using *CE property 7*, Eq. (3.1.14b), we have ([4], §9.1):

$$P(\Omega \mid Y_n(X_1,...,X_n) \mid X_1,...,X_n) = Y_n(X_1,...,X_n), \text{ or } P(\Omega \mid Y_n \mid X_1,...,X_n) = Y_n \tag{4.1.4}$$

Furthermore, using *CE property 5*, Eq. (3.1.12), we have:

$$\langle Y_{n+1}\rangle \underset{Eq.(2.1.12)}{=} \langle P(\Omega \mid Y_{n+1} \mid X_1,...,X_n)\rangle \underset{Eq.(3.1.2)}{=} \langle Y_n\rangle$$

Therefore, we have the following property for *M.G.* $Y_n$ ([4], §9.1):

$$\langle Y_n\rangle = \langle Y_0\rangle, \quad \forall n \ge 0 \tag{4.1.5}$$

**Example 4.1.1** ([4], example 9.1.3): Assume that $\{X_n, n \ge 0\}$ is a homogeneous Markov chain (MC), its transition matrix is $P = (p_{ij})$; Assume also that $\{\varphi(i)\}$ is a non-negative bounded sequence, satisfying:

$$\varphi(i) = \sum_j p_{ij}\, \varphi(j) \underset{(2.4.18b)}{=} \sum_j p(j \mid i)\, \varphi(j) \tag{4.1.6}$$

Let $Y_n = \varphi(X_n)$, then $\langle |Y_n|\rangle < +\infty$ (because $\{\varphi(i)\}$ is finite). Also:

$$P(\Omega \mid Y_{n+1} \mid X_0,\dots X_n) \underset{by\,def.}{=} P(\Omega \mid \varphi(X_{n+1}) \mid X_0,\dots X_n) \underset{Markov+(2.1.5a)}{=} P(\Omega \mid \varphi(X_{n+1}) \mid X_n) \tag{4.1.7a}$$

On the other hand:

$$\begin{aligned} P(\Omega \mid \varphi(X_{n+1}) \mid X_n = i) &= \sum_j P(\Omega \mid \varphi(X_{n+1}) \mid X_{n+1} = j) \cdot P(X_{n+1} = j \mid X_n = i) \\ &\underset{(2.4.10a)}{=} \sum_j \varphi(j)\, p(j \mid i) \underset{(4.1.6)}{=} \varphi(i) \end{aligned} \tag{4.1.7b}$$

$$\begin{aligned} &\therefore P(\Omega \mid \varphi(X_{n+1}) \mid X_n) = \varphi(X_n) \underset{\text{By Defi.}}{=} Y_n \\ &\therefore P(\Omega \mid Y_{n+1} \mid X_0,\dots X_n) \underset{(4.1.7a)}{=} P(\Omega \mid \varphi(X_{n+1}) \mid X_n) = Y_n \quad \text{QED} \end{aligned} \tag{4.1.8}$$

Therefore, $\{Y_n\}$ is a martingale about $\{X_n\}$.

**Example 4.1.2** (see [4], example 9.1.4): Assume that $\{X_n, n \ge 0\}$ is a homogeneous Markov chain (MC), its transition matrix is $P = (p_{ij})$; Assume that λ is a non-zero eigenvalue of matrix *P*, corresponding to eigenvector $\{\varphi(i)\}$, as in Eq. (2.4.18b), i.e.

$$\lambda\varphi(i) = \sum_j p_{ij}\, \varphi(j) = \sum_j p(j \mid i)\, \varphi(j), \quad \forall i \tag{4.1.9}$$

Assume also that $\langle |\varphi(X_n)|\rangle < +\infty, \forall n > 0$,

and let: $Y_n = \lambda^{-n}\varphi(X_n)$, then $\{Y_n\}$ is a *M.G* about $\{X_n\}$ (4.1.10)

*Proof:* By definition and Markov property, we have:

$$P(\Omega \mid Y_{n+1} \mid X_0, \dots X_n) \underset{(4.1.10)}{=} P(\Omega \mid \lambda^{-n-1}\varphi(X_{n+1}) \mid X_0, \dots X_n) \underset{Markov}{=} \lambda^{-n-1} P(\Omega \mid \varphi(X_{n+1}) \mid X_n)$$

But, using unit operator and definition of Markov matrix, we have:

$$P(\Omega \mid \varphi(X_{n+1}) \mid X_n = i) = \sum_j P(\Omega \mid \varphi(X_{n+1}) \mid X_{n+1} = j) P(X_{n+1} = j \mid X_n = i)$$

$$\underset{(2.4.10a)}{=} \sum_j \varphi(j)\, p(j \mid i) = \sum_j p_{ij}\, \varphi(j) \underset{(4.1.9)}{=} \lambda \varphi(i)$$

Therefore, $P(\Omega \mid \varphi(X_{n+1}) \mid X_n) = \lambda \varphi(X_n)$, and we have the condition in (4.1.2):

$$P(\Omega \mid Y_{n+1} \mid X_0, \dots X_n) = \lambda^{-n-1} P(\Omega \mid \varphi(X_{n+1}) \mid X_n) = \lambda^{-n} \varphi(X_n) \underset{(4.1.10)}{\equiv} Y_n \qquad \text{QED}$$

This example can be easily extended to time-discrete, state-continuous and homogeneous MC (see Eq. (2.4.19)).

**Example 4.1.3** (***Doob Martingale*** see [4], example 9.1.7): Assume that $\{X_n, n \ge 0\}$ is a sequence of $R.V$, $Y$ is a $S.V$ and $\langle | Y | \rangle < +\infty$, then the sequence of $R.V$ defined by: $Y_n = P(\Omega \mid Y \mid X_0, X_1, \cdots, X_n)$ is a martingale about $\{X_n\}$, which is called a Doob martingale.

*Proof*: $\langle | Y_n | \rangle = \langle | P(\Omega \mid Y \mid X_0, X_1, \cdots, X_n) | \rangle \le \langle P(\Omega \mid \{| Y |\} \mid X_0, X_1, \cdots, X_n) \rangle \underset{(3.1.12b)}{=} \langle | Y | \rangle$

*Moreover*, $\mathrm{P}(\Omega \mid \{Y_{n+1}\} \mid X_0, X_1, \cdots, X_n)$

$$\underset{By\,def.}{=} P(\Omega \mid \{P(\Omega \mid Y \mid X_0, X_1, \cdots, X_n)\} \mid X_0, X_1, \cdots, X_{n+1}) \underset{(3.1.18)}{=} P(\Omega \mid Y \mid X_0, X_1, \cdots, X_n) = Y_n$$

**Example 4.1.4** (***Wald Martingale*** see [4], example 9.1.8): Assume that $\{X_0 = 0, X_1, X_2, \cdots\}$ is a mutual independent, identically distributed (i.i.d) sequence of $R.V$; assume also that there exists a finite-moment generating function: $\varphi(\lambda) = \langle e^{\lambda Y_1} \rangle$ for some $\lambda \ne 0$. Define $Y_0 = 1, Y_n = [\varphi(\lambda)]^{-n} \cdot \exp\{\lambda(X_1 + \cdots + X_n)\}, n \ge 1$; then $\{Y_n\}$ is a martingale about $\{X_n\}$.

*Proof*: First, we verify Eq. (4.1.2a):

$$\langle | Y_n | \rangle = [\varphi(\lambda)]^{-n} \cdot \langle \exp\{\lambda(X_1 + \cdots + X_n)\} \rangle \underset{X_i\ independent}{=} [\varphi(\lambda)]^{-n} \cdot \sum_{i=1}^{n} \langle \exp\{\lambda X_i\} \rangle$$

$$\underset{identical\ distribution}{=} [\varphi(\lambda)]^{-n} \cdot [\exp\{\lambda X_1\}]^n \underset{by\ definition}{=} 1 < +\infty$$

Next, we verify Eq. (4.1.2b):

$$P(\Omega \mid \{Y_{n+1}\} \mid X_0, X_1, \cdots, X_n) \underset{by\,def.}{=} [\varphi(\lambda)]^{-n-1} \cdot P(\Omega \mid \exp\{\lambda(X_1 + \cdots + X_{n+1})\} \mid X_0, X_1, \cdots, X_n)$$

$$\underset{(3.1.13)}{=} [\varphi(\lambda)]^{-n-1} \cdot \exp\{\lambda(X_1 + \cdots + X_n)\} \cdot P(\Omega \mid \exp\{\lambda(X_{n+1})\} \mid X_0, X_1, \cdots, X_n)$$

$$\underset{X_i\ independent}{=} [\varphi(\lambda)]^{-n-1} \cdot \exp\{\lambda(X_1 + \cdots + X_n)\} P(\Omega \mid \exp\{\lambda(X_{n+1}) \mid \Omega)$$

$$\underset{by\,def.}{\equiv} Y_n[\varphi(\lambda)]^{-1} \langle \exp\{\lambda(X_{n+1})\rangle \underset{ident.\,dist.}{=} Y_n[\varphi(\lambda)]^{-1} \langle \exp\{\lambda(X_1)\rangle = Y_n[\varphi(\lambda)]^{-1} \varphi(\lambda) = Y_n$$

### 4.2. Discrete-time Martingale Based on Filtration

A *martingale* (*M.G.)* based on filtration consists of a filtration and an adapted *S.P*. which has the property of being a "fair game", that is, the expected future reward given current information is exactly the current value of the process. We now discuss related definitions, propositions and examples using *PBN*.

**Definition 4.2.1** ([5], §4.1): A ***martingale*** is a pair $(X_n, \Im_n)$, where $\{\Im_n\}$ is a filtration and $\{X_n\}$ an *integrable* (i.e. $\langle | X_n | \rangle < \infty$) *S.P.* adapted to this filtration such that

$$E(X_{n+1} \mid \Im_n) \equiv P(\Omega \mid X_{n+1} \mid \Im_n) = X_n \tag{4.2.1}$$

**Proposition 4.2.1**: If $X_n = \sum_{i=0}^{n} D_i$ then the canonical filtration for $\{X_n\}$ is the same as the canonical filtration for $\{D_n\}$. Further, $\{X_n, \Im_n\}$ is a martingale if and only if $\{D_n\}$ is an integrable *S.P.*, adapted to $\{\Im_n\}$, such that $P(\Omega \mid D_{n+1} \mid \Im_n) = 0$ almost sure for all *n*.

*Proof*: Since the transformation from $(X_0, ..., X_n)$ to $(D_0, ..., D_n)$ is continuous and invertible, it follows (from [5], Corollary 1.2.17) that for each *n*. By Definition 2.3.3 we see that $\{X_n\}$ is adapted to a filtration $\{\Im_n\}$ if and only if $\{D_n\}$ is adapted to this filtration. It is very easy to show by induction on *n* that $\langle | X_k | \rangle < \infty$ for $k = 0, \ldots, n$ if and only if $\langle | D_k | \rangle < \infty$ for $k = 0, \ldots, n$. Hence, $\{X_n\}$ is an integrable *S.P.* if and only if $\{D_n\}$ is. Finally, with $X_n$ measurable on $\Im_n$ it follows from the linearity of the *C.E.* that

$$P(\Omega \mid X_{n+1} \mid \Im_n) - X_n \underset{(3.2.6)}{=} P(\Omega \mid X_{n+1} - X_n \mid \Im_n) \underset{\text{by def.}}{=} P(\Omega \mid D_{n+1} \mid \Im_n) \tag{4.2.2}$$

And the alternative expression for the martingale property follows from Eq. (4.2.1).

**Example 4.2.1**: The random walk $S_n = \sum_{k=1}^{n} \xi_k$ with $\xi_k$ independent, identically distributed (i.i.d), such that $\langle | \xi_k | \rangle < \infty$ and $\langle \xi_k \rangle = 0$, is a MG (for its canonical filtration). More generally, $\{S_n\}$ is a MG even when the independent and integrable *R.V.* $\xi_k$ of zero mean have non-identical distributions. Further, the canonical filtration may be replaced by the filtration $\sigma(\xi_1,\ldots,\xi_n)$. Indeed, this is just an application of *Proposition 4.2.1* for the case where the differences $D_k = S_k - S_{k-1} = \xi_k,\ k \ge 1$ (and $D_0 = 0$), are independent, integrable and $P(\Omega \mid D_{n+1} \mid D_0, D_1, \ldots D_n) = P(\Omega \mid D_{n+1} \mid \Omega) = 0$ for all $n \ge 0$ by our assumption that $\xi_k$ is independent and $\langle \xi_k \rangle = 0$ for all $k$.

Alternatively, for a direct proof recall that $\langle | S_n | \rangle \le \sum_{k=1}^{n} \langle | \xi_n | \rangle < \infty$ for all $n$, that is, the *S.P* $S_n$ is integrable. Moreover, since $\xi_{n+1}$ is independent of $\sigma(S_1,\ldots,S_n)$, $\langle \xi_{n+1} \rangle = 0$, and $P(\Omega \mid S_n \mid S_1,\ldots,S_n) \underset{Eq.(3.1.14b)}{=} S_n$, we have that

$$P(\Omega \mid S_{n+1} \mid S_1,\ldots S_n) = P(\Omega \mid S_n \mid S_1,\ldots S_n) + P(\Omega \mid \xi_{n+1} \mid S_1,\ldots S_n) = S_n + \langle \xi_{n+1} \rangle = S_n \quad (4.2.3)$$

This is implying that $\{S_n\}$ is a *MG* for its canonical filtration.

**Example 4.2.2**: Let $Y_n = \sum_{k=1}^{n} V_k \xi_k$, where $\{V_n\}$ is an arbitrary bounded *S.P.* and $\{\xi_n\}$ is a sequence of integrable *R.V.* such that for $n$ = 0, 1, …, both $\langle \xi_{n+1} \rangle = 0$ and $\xi_{n+1}$ is independent of $\Im_n = \sigma(\xi_1,\ldots,\xi_n,V_1,\ldots,V_n)$. Then, $\{Y_n\}$ is a *MG* for its canonical filtration and even for the possibly larger filtration $\{\Im_n\}$. This is yet another application of Proposition 4.1.7, now with the differences $D_k = Y_k - Y_{k-1} = V_k \xi_k,\ k \ge 1$ (and $D_0 = 0$). Indeed, we assumed $\xi_k$ are integrable and $|V_k| \le C_k$ for some non-random finite constants $C_k$, resulting with $\langle | D_k | \rangle \le C_k \langle | \xi_k | \rangle < \infty$, whereas:

$$\begin{aligned} P(\Omega \mid D_{n+1} \mid \Im_n) &\underset{By\ Definition}{=} P(\Omega \mid V_{n+1}\xi_{n+1} \mid \Im_n) \underset{Eq.(3.2.4)}{=} V_{n+1} P(\Omega \mid \xi_{n+1} \mid \Im_n) \\ &\underset{Eq.(3.2.9)}{=} V_{n+1} \langle \xi_{n+1} \rangle \ (\because \xi_{n+1} \text{ is independent of } \Im_n) = 0 \text{ (by given condition)} \end{aligned} \quad (4.2.4)$$

This gives us the martingale property.

A special case of *Example* 4.2.2 is when the auxiliary sequence $\{\xi_k\}$ is independent of the given *S.P.* $\{V_n\}$ and consists of zero-mean, independent, identically distributed (i.i.d) *R.V.* For example, random i.i.d. signs $\xi_k \in \{1,-1\}$ (with $P(\xi_k = 1) = P(\xi_k = 1 \mid \Omega) = \frac{1}{2}$) are commonly used in discrete mathematics applications.

## 4.3. Continuous-Time Martingale Based on Filtration

**Definition 4.3.1** (see [8]): A **continuous-time *M.G.* with respect to** the stochastic process $X_t$ is a stochastic process $Y_t$ such that for all $t$:

$$E(|Y_t|) \equiv \langle |Y_t| \rangle \equiv P(\Omega | \{|Y_t|\} | \Omega) < \infty \quad (4.3.1a)$$

$$E(Y_t | \{X_\tau, \tau \le s\}) \equiv P(\Omega | Y_t | \{X_\tau, \tau \le s\}) = Y_s, \ a.s. \text{ for all } s \le t \quad (4.3.1b)$$

This expresses the property that the conditional expectation of an observation at time *t*, given all the observations up to time *s*, is equal to the observation at time *s* (of course, provided that $s \le t$). If we use our notation of filtration $\mathsf{F}$ in §2.3 (or [6], §3.1), we have:

**Definition 4.3.2** (see [6], §3.1): Let $(\Omega, \Im, P)$ be a probability space, $\mathsf{F} = \{\Im_t, t \ge 0\}$ a filtration contained in $\Im$; and $\{Z_t \equiv Z(t), t \ge 0\}$ a *stochastic process adapt*ed to $\mathsf{F}$ ($Z_t$ *is* $\Im_t$*-measurable* for all *t* ). Then $\{[Z_t, \Im_t], t \ge 0\}$ is *a martingale* (or $\{Z_t\}$ is a martingale with respect to $\mathsf{F}$ ) if for all $t \ge 0$,

$$E(|Z_t|) \equiv \langle |Z_t| \rangle \equiv P(\Omega | \{|Z_t|\} | \Omega) < \infty \quad (4.3.2a)$$

$$E(Z_t | \Im_s) \equiv P(\Omega | Z_t | \Im_s) = Z_s, \text{ a.s. for all } 0 \le s < t \quad (4.3.2b)$$

Note that the time index t may be discrete or continuous and the horizon may be finite or infinite. Also note that the definition involves the filtration $\mathsf{F} = \{\Im_t\}$ in an integral way. If the filtration is understood, however, we may say simply that $\{Z_t\}$ is a martingale.

**Definition 4.3.3** (see [7], §4.6): If, instead of (4.3.2b),
$P(\Omega | Z_t | \Im_s) > Z_s$, for all $0 \le s < t$ a.s. – then $\{Z_t\}$ is a *submartingale* with respect to $\{\Im_t\}$
$P(\Omega | Z_t | \Im_s) \le Z_s$, for all $0 \le s < t$ a.s. – then $\{Z_t\}$ is a *supermartingale* with respect to $\{\Im_t\}$

The first and most obvious way to construct a martingale is as the sum of independent mean-zero random variables. Indeed, the idea originated as the stochastic process describing the wealth of a gambler playing games of chance at fair odds.

**Example 4.3.1** (***Gambler's wealth***; see [6], §3.1, Example A):
Let $\{X_1, \ldots X_i, \ldots\}$ be a (finite or infinite) sequence of independent random variables on the probability space $(\Omega, \Im, P)$, each with mean zero and with finite absolute deviation: $\langle X_i \rangle = 0$ and $\langle |X_i| \rangle < +\infty$, all *i*. Let $Z_0 = 0$, and for each *k* = 1, 2, ..., define the random variable $Z_k = \sum_{i=1}^k X_i$ to be the partial sum of the first k elements in the sequence. For each *k* let $\Im_k \subseteq \Im$ be the smallest σ-algebra for which $\{Z_j\}_{j=1}^k$ are measurable, and let $\mathsf{F} = \{\Im_k\}$ be the filtration consisting of this (increasing) sequence. Then the stochastic process $\{Z_k\}$ is a martingale on the filtered space $(\Omega, \mathsf{F}, P)$.

$$Proof: \langle |Z_k| \rangle = \langle \{|\textstyle\sum_{i=1}^k X_i|\} \rangle \le \langle \textstyle\sum_{i=1}^k |X_i| \rangle = \textstyle\sum_{i=1}^k \langle |X_i| \rangle \underset{given}{<} +\infty \quad (4.3.3a)$$

Moreover, note that for any $j < k$, $E(Z_k | \Im_j) \equiv P(\Omega | Z_k | \Im_j) \underset{By\,def.}{=} P(\Omega | [\sum_{i=1}^k X_i] | \Im_j)$

$$\underset{By\,def.}{=} P(\Omega|[\sum_{i=1}^{j} X_i]|\Im_j) + P(\Omega|[\sum_{i=j+1}^{k} X_i]|\Im_j) \underset{(3.2.11)}{=} \sum_{i=1}^{j} X_i + \sum_{i=j+1}^{k} P(\Omega|X_i|\Im_j)$$

$$\underset{(3.2.12)}{=} \sum_{i=1}^{j} X_i + \sum_{i=j+1}^{k} \langle X_i \rangle \underset{\langle X_i \rangle = 0}{=} \sum_{i=1}^{j} X_i \tag{4.3.3b}$$

***Poisson Process*** (see [2], §5.1) $N(t) \equiv N_t$ with intensity λ has following properties:
(1). $\{N(t), t \ge 0\}$ is non-negative process with independent increments and $N(0) = 0$;
(2). It is homogeneous and its probability distribution is given by:

$$m(k,t) \equiv P(N_t = k|\Omega) = P(k|\Omega(t)) \underset{\substack{Poisson\\Distribution}}{=} \frac{(\lambda t)^k}{k!} e^{-\lambda t} \tag{4.3.4a}$$

Using identity operator, one can easily find that:

$$\mu(t) \equiv \langle N(t) \rangle = \lambda t; \quad \sigma^2(t) \equiv P(\Omega|[N(t) - \bar{N}(t)]^2|\Omega) = \lambda t \tag{4.3.4b}$$

$$[\lambda]_{P.D} = [t]_{P.D}^{-1} = T^{-1}, \quad [N_t]_{P.D} = [|k,t)]_{P.D} = [m(k,t)]_{P.D} = 1 \tag{4.3.4c}$$

**Example 4.3.2:** Let $N_t$ be a Poisson process with intensity $\lambda$ and $X_t = N_t - \lambda t$, then $X_t$ is a martingale.

$$Proof: P(\Omega|X_t|\Im_s) = P(\Omega|\{N_s + N_t - N_s\}|\Im_s) - \lambda t \underset{(3.2.11)}{=} N_s + P(\Omega|\{N_t - N_s\}|\Im_s) - \lambda t$$

$$\underset{(3.2.13)}{=} N_s + P(\Omega|\{N_t - N_s\}|\Omega) - \lambda t = N_s - \lambda s \equiv X_s \tag{4.3.5}$$

**Lemma 4.3.1** (see [7], §5.1): More generally, if $X_t$ is *any bounded S.P with independent increments*, then an identical calculation shows that $Z_t = X_t - \langle X_t \rangle$ is a martingale.

$$P(\Omega|Z_t|\Im_s) = P(\Omega|\{X_s + X_t - X_s\}|\Im_s) - \langle X_t \rangle \underset{(3.2.11)}{=} X_s + P(\Omega|\{X_t - X_s\}|\Im_s) - \langle X_t \rangle$$

$$\underset{(3.2.13)}{=} X_s + P(\Omega|\{X_t - X_s\}|\Omega) - \langle X_t \rangle = X_s - \langle X_s \rangle \equiv Z_s \tag{4.3.6}$$

**Wiener Process** (see [1], §5.1; [6], §1.1; [10], §10.1)**:** It is a homogeneous process $\{W(t), t \ge 0\}$ with independent increments and $W(0) = 0$. Its probability density is a normal distribution $N(0, t\sigma^2)$:

$$f(x,t) \equiv P(x;t|\Omega) \equiv P(x|\Omega(t)) \underset{\substack{Normal\\Distribution}}{=} \frac{1}{\sqrt{2\pi t}\,\sigma} \exp[-\frac{x^2}{2t\sigma^2}] \tag{4.3.7}$$

$$\tilde{\mu}(t) \equiv P(\Omega|W(t)|\Omega) \equiv \langle W_t \rangle = 0, \quad \tilde{\sigma}^2(t) \equiv P(\Omega|W(t)^2|\Omega) \equiv \langle W_t^2 \rangle = t\sigma^2 \tag{4.3.8}$$

Because Wiener process has independent increments and has zero mean, it is a martingale according to Lemma 4.3.1.

**Brownian motion** (see [1], §5.1; [6] §1.1; [7], Brownian motion; [10], §10.3)**:** It is associated with a *standard Winner process* $W_s(t)$ (with σ = 1) as follows:

$$X(t) = \mu t + \sigma W_s(t) \tag{4.3.9}$$

Using Eq. (4.3.8) (with σ = 1), one can easily find that:

$$Drift:\quad \tilde{\mu}(t) \equiv \langle X_t \rangle = \mu t; \quad Variance:\quad \tilde{\sigma}^2(t) \equiv \langle [X_t - \tilde{\mu}(t)]^2 \rangle = t\sigma^2 \tag{4.3.10}$$

**Example 4.3.3:** Because Brownian motion has independent time increments, according to Lemma 4.3.1, we have the following martingale $Z(t)$:

$$Z(t) \equiv X(t) - \langle X(t) \rangle = X(t) - \mu t, \quad E(Z_t \,|\, \Im_s) \equiv P(\Omega \,|\, Z_t \,|\, \Im_s) = Z_s \tag{4.3.11}$$

$$\langle Z(t) \rangle = 0 \quad (all\ t \ge 0), \quad \langle (Z_t - Z_s)^2 \rangle = \sigma^2 (t-s) \quad (all\ t > s \ge 0) \tag{4.3.12}$$

**Example 4.3.4:** Let $X_t$ be a Brownian motion given by Eq. (4.3.9), then we have a martingale: $M_t \equiv (X_t - \mu t)^2 - \sigma^2 t \underset{(4.3.11)}{=} Z_t^2 - \sigma^2 t$ (4.3.13)

$$Proof:\quad Let\ \Delta Z_{t-s} \equiv Z_t - Z_s\ for\ t > s, \quad then: P(\Omega \,|\, \Delta Z_{t-s} \,|\, \Im_s) = \langle \Delta Z_{t-s} \rangle \underset{(4.3.12)}{=} 0 \tag{4.3.14a}$$

$$Also:\quad P(\Omega \,|\, \Delta Z_{t-s}{}^2 \,|\, \Im_s) \underset{independence}{=} \langle (Z_t - Z_s)^2 \rangle \underset{(4.3.12)}{=} \sigma^2 (t-s) \tag{4.3.14b}$$

$$Now:\quad E(M_t \,|\, \Im_s) \equiv P(\Omega \,|\, M_t \,|\, \Im_s) = P(\Omega \,|\, \{Z_s + \Delta Z_{t-s}\}^2 \,|\, \Im_s) - \sigma^2 t$$

$$= P(\Omega \,|\, Z_s{}^2 \,|\, \Im_s) + P(\Omega \,|\, \Delta Z_{t-s}{}^2 \,|\, \Im_s) - 2P(\Omega \,|\, Z_s \Delta Z_{t-s} \,|\, \Im_s) - \sigma^2 t$$

$$\underset{\substack{(4.3.14b)\\(3.2.11)}}{=} Z_s{}^2 + \sigma^2 (t-s) - 2 Z_s \langle \Delta Z_{t-s} \rangle - \sigma^2 t \underset{(4.3.12)}{=} Z_s{}^2 - \sigma^2 s \underset{(4.3.13)}{=} M_s \tag{4.3.15}$$

***Dimensional Analysis***: Assuming Wiener process $[W]_{P.D} = L$ and $[t]_{P.D} = T$, we see that Brownian process $[X]_{P.D} = L$ and we have the following *P.D* for related quantities:

$$[f(x;t)]_{P.D} = L^{-1}, \quad [|\,\Omega_t)] = [|\,\Im_t)]_{P.D} = L^{-1/2}, \quad [(\Omega\,|]_{P.D} = L^{1/2} \tag{4.3.16a}$$

$$[\mu]_{P.D} = LT^{-1}, \quad [\sigma] = LT^{-1/2}, \quad [Z_t]_{P.D} = L, \quad [M_t]_{P.D} = L^2 \tag{4.3.16b}$$

Note: *standard Wiener process* $W_s(t)$ is defined by setting in Eq. (4.3.7) with σ = 1. This is just a selection of unit: the position here is measured by using a quantity $X_s = X/\sigma$, where σ is given by Eq. (4.3.7). In this case, we have $[W_s]_{P.D} = [W/\sigma]_{P.D} = T^{1/2}$, which is consistent with Eq. (4.3.9) and Eq. (4.3.16): $[X]_{P.D} = [\mu t]_{P.D} = [\sigma W_s]_{P.D} = L$.

## ***Summary***

In this paper, we continued to investigate *PBN* (*Probability Bracket Notation*) in probability modeling. We discussed probability space, σ-fields, filtration and the symbols in *PBN* with probability spaces associated with *R.V* and *S.P*. Then we studied important properties of conditional expectation with respect to *R.V* or to filtration in *PBN*. Finally, we explored the definitions and examples of discrete or continuous martingales based on *R.V* or based on filtration in *PBN*. We saw that, in each topic, *PBN* simplified the expression and manipulation of related formulas. We also applied physical dimensional analysis to verify our propositions in *PBN*.

*PBN* has adapted many concepts and symbols from the Dirac notation in QM. We demonstrated that the Dirac delta function, the unit operator and the indicator operator make it possible to investigate certain probability issues without explicit usage of Measure theory.

Moreover, our recent study [11] reveals that many expressions in the Dirac notation are naturally shifted to expressions in *PBN* under Wick rotations. Therefore, the two notations may have deeper co-relation and might be the two faces of the same coin.